\documentclass[12pt]{article}
\usepackage{epsfig}
\usepackage{rotating}
\title{Patterns of Geodesics, Shearing, and Anosov Representations of
  the Modular Group}
\author{Richard Evan Schwartz \thanks{\hskip 5 pt 
Supported by 
N.S.F. Grant DMS-2505281}}

\newtheorem{theorem}{Theorem}[section]

\newtheorem{lemma}[theorem]{Lemma}

\newtheorem{corollary}[theorem]{Corollary}

\def\startproof{{\bf {\medskip}{\noindent}Proof: }}

\def\endproof{$\spadesuit$  \newline}

\def\H{\mbox{\boldmath{$H$}}}% 
\def\P{\mbox{\boldmath{$P$}}}% 
\def\R{\mbox{\boldmath{$R$}}}% 
\def\Z{\mbox{\boldmath{$Z$}}}% 

\begin{document}

\maketitle

\begin{abstract}
  Let $X=SL_3(\R)/SO(3)$.
  Let $\cal DFR$ be the
 space of discrete 
  faithful representations of the
  modular group into
  ${\rm Isom\/}(X)$ which map the order
  $2$ generator to an isometry with a unique fixed point.
  This paper continues my work on the
  Barbot component $\cal B$ of $\cal DFR$,
  showing how all the representations in
  $\cal B$ are shears of the Pappus modular
  groups in an appropriate sense.
\end{abstract}

      \section{Introduction}

Let  $X=SL_3(\R)/SO(3)$.  This is a prototypical
higher rank symmetric space.
In this paper we consider the moduli space
$\mathcal DFR$ of conjugacy classes of
discrete and faithful
representations of the modular group
$\Z/3* \Z/2$ into
${\rm Isom\/}(X)$ which map
the order $2$ elements to isometries
having a unique fixed point in $X$.

The {\it Pappus representations\/}
  are a $2$-parameter subfamily
  of $\mathcal DFR$ which I constructed in
  my 1993 paper
  [{\bf S0\/}] and then revisited
  in my recent paper [{\bf S1\/}].
  These groups exhibit many features
  that, much later and more generally,
  appeared in
   higher Teichmuller Theory, e.g. in [{\bf Lab\/}],
  [{\bf GW\/}]. [{\bf Bar\/}],
  [{\bf BCLS\/}], and [{\bf KL\/}].

The Pappus representations 
  arise as the projective symmetry groups
  of {\it convex marked box orbits\/}.
  A {\it convex marked box\/} is convex quadrilateral
  with two additional marked points on a pair
  of opposite edges. One gets a marked box
  orbit by starting with one
    marked box and iteratively applying
  operations that are derived from
  Pappus's Theorem, a classic theorem
  in projective geometry.
  
The Pappus representations
are
  nowadays classified as {\it relatively Anosov groups in the
    Barbot component\/}.
  This point is view is exposited in
  [{\bf BLV\/}] and [{\bf KL\/}].
  Let ${\mathcal P\/} \subset \mathcal DFR$ denote
  the subset consisting of Pappus modular
  group representations.  The component
  $\mathcal B$ of $\mathcal DFR$ containing
  $\mathcal P$ is called the {\it Barbot component\/}.
  It is partially understood
  thanks to  [{\bf S0\/}] and [{\bf BLV\/}].
  
  In [{\bf BLV\/}], T. Barbot,
  G.-S. Lee, and V. P. Valerio build
    on [{\bf S0\/}] and construct
  a $3$-parameter family of Anosov
  representations which are defined
  in terms of modified operations on
  marked boxes.  Using their
  {\it morphed marked boxes\/} (my terminology)
  they construct a $4$-parameter family of
  representations of $\Z/3 * \Z/3$ into
  $SL_3(\R)$, all of which are Anosov.
  Using an Implicit Function argument,
  they show that a $3$-parameter subset
  of these are index $2$ subgroups of Anosov
  elements of $\cal DFR$.  I gave \footnote{The original version
    of this paper contained this analysis. However, I split the
    original version of this paper into two halves, offloading
    the analysis of the Barbot-Lee-Valerio construction to a
    separate paper. I did this because the original paper was
    very long.}
  a complete
  algebraic analysis of their construction
  in [{\bf S2\/}].
  
  Now I will summarize the results of
  [{\bf S2\/}].  Let $\cal R$ be the space
  of representations of $\Z/3*\Z/2$
  with the following property:
  The order $2$ generator is mapped to an isometry
  with a unique fixed point in $X$, and this point
      is not contained in the fixed set of the order $3$ generator.
  
\begin{theorem}
  \label{one}
  The space $\R$ is homeomorphic to $\R^3-\{0\}$.
  The set $\cal P$ of Pappus representations is a properly embedded
  plane in $\cal R$.  The connected component of
  ${\cal R\/}-{\cal P\/}$ that is disjoint from the origin
  consists entirely of Anosov representations in $\cal DFR$.  The
  union of $\cal P$ with this Anosov component is
  a connected component of $\cal DFR$.
\end{theorem}

Let ${\cal B\/} \subset \cal R$ be the union of
$\cal P$ and the Anosov component of ${\cal R\/}-{\cal P\/}$.
In this paper, I will give alternate descriptions of the groups
in $\cal B$ as isometry groups of certain patterns of geodesics
and flats in $X$.  There are two descriptions
per representation, distinct in the generic case.
I call these the
{\it prism descriptions\/}.  (The end of this introduction
explains the name.)
The prism descriptions and the marked box description
complement each other:  Properties of the parametrization,
such as injectivity, are easy to see from the prism
descriptions, and the discreteness/Anosov properties of
the individual representations are
easier to see from the morphed marked box description.
In particular, all the hard work concerning the Anosov
nature of the representations is already done in [{\bf BLV\/}].

The prism descriptions are an elaboration of my paper
  [{\bf S0\/}], in which  I re-interpreted
  the Pappus modular groups
  as isometry groups of what I
  call {\it Farey patterns\/}.
  These are embedded patterns
  of geodesics in $X$ that have the
  same asymptotic structure as the
  Farey triangulation.  Each geodesic in the
  Farey pattern is a
  {\it medial geodesic\/}.   These are
  geodesics which are
    angle bisectors of the Weyl chambers
    within the flats that contain them.
    These patterns do not generally lie in a
    totally geodesic slice of $X$.  They are
    bent, sort of like pleated planes.

    \begin{theorem}
      \label{extra}
      Let $\rho \in {\mathcal B\/}-{\cal P\/}$ be arbitrary.
      Let $\Gamma=\rho(\Z/3*\Z/2)$.
      Then there are two infinite embedded patterns
      $Y_{\rho,1}$ and $Y_{\rho,2}$ of
      geodesics in $X$ such that
      $\Gamma \subset {\rm Isom\/}(Y_{\rho,1})$
      and $\Gamma \subset {\rm Isom\/}(Y_{\rho,1})$
      with finite index.
      In the generic case, $\Gamma$ equals these
      isometry groups, and
      $Y_{\rho_1}, Y_{\rho,2}$ are not isometric to each other.
      Finally, the pattern of flats containing
      $Y_{\rho,k}$ is embedded for $k=1,2$.
    \end{theorem}

        Just as in the classic Farey triangulation, the geodesics
    in the Farey patterns are organized into
        triples of mutually asymptotic geodesics which
        I call {\it triangles\/}.  (The triple of flats containing the
        geodesics in a triangle are the {\it prisms\/}.) Two triangles
        $\tau_1,\tau_2$ are
        {\it adjacent\/} if they have a geodesic
        $\gamma$ in common.  The pair
        $(\tau_1,\tau_2')$ is a {\it shear\/}
        of $(\tau_1,\tau_2)$ if $\tau_2'=I(\tau_2)$
        where $I$ is an isometry of $X$ which translates
        along $\gamma$.  There is a
        $1$-parameter family of such translations.
        The amount of shearing, which we call
        the {\it strength\/} of the shear, is the translation
        length of $I$.
        More generally, if $J$ is an isometry of $X$,
        the pair $(J(\tau_1),J(\tau_2'))$ is
        a  shear of $(\tau_1,\tau_2)$.

        Speaking more globally, 
        a {\it shear\/}
        of a Farey pattern is another pattern of geodesics
        in which adjacent triangles in the new pattern are
        shears of adjacent triangles in the original
        pattern -- all of the same strength and in
        the same direction so to speak.

    \begin{theorem}
      \label{three}
      The space $\mathcal B$ has two different $1$-dimensional
      foliations by proper and canonically parametrized
      rays. Each ray has its endpoint
      in $\mathcal P$.  A representation $\rho \in \mathcal B$
      lies $d$ units along a ray in the first foliation
        if and only if it lies $d$ units
      along a ray in the second foliation.
      In this case, the patterns
      $Y_{\rho,1}$ and $Y_{\rho,2}$ are each
      strength-$d$ shears respectively of
      Farey patterns $F_{\rho,1}$ and
      $F_{\rho,1}$.  Generically, $F_{\rho,1}$ and
      $F_{\rho,2}$ are not isometric to either.
    \end{theorem}

        For each $d$, the double foliation sets up a map
        $\phi_d: {\mathcal P\/} \to {\mathcal P\/}$, which I
        call the {\it shearing dynamics\/}.  Here is a description.
        We start with a Pappus
        representation $p_0$, and then move $d$-units along the
        ray in the first foliation that contains it.   We arrive at
        some representation $p_d$. Then we move
        $d$ units backwards along the ray in the second foliation
        that contains $p_d$. This gives us $\phi_d(p)$.
        The map $\phi_d$ is
        a self-homeomorphism of $\mathcal P$ that is
        the product of $2$ involutions.  It is
        the {\it bounce map\/} associated to the initial
        segments of length $d$ of our two ray foliations.
        
        I did a little bit of experimentation with these shearing
        dymamics.  (I am not yet completely confident in my
        way of computing it.)  For the parameters I 
        looked it, it seems that
        $\phi_d$ has infinite order and that its orbits
        are dense in finite unions of $1$-dimensional curves.
                \newline
    
    This paper is organized as follows.
    \begin{itemize}
\item   In \S 2, I discuss classic shearing
  of the modular group in the hyperbolic plane.
  The shearing phenomenon we uncover in
  $X$ extends what happens in $\H^2$.
  
\item  In \S 3, I give some background material
  about the space $X$.
  
\item  In \S 4, I  explore the geometry of {\it prisms\/},
  namely triples of flats which contain mutually
  asymptotic medial geodesics.  (I only consider
  prisms defined by flats having negative triple
  invariants; the other kinds of prisms are presumably
  related to the Goldman-Hitchin component
  [{\bf G1\/}], [{\bf Hit\/}]
  of $\mathcal DFR$.)
  At the end of \S 4, I describe a space
  ${\mathcal A\/}$ of pairs $(\Pi,p)$ where $\Pi$ is a prism
  and $p \in \Pi$ is some point.
  There is a map $\rho: {\mathcal A\/} \to {\mathcal R\/}$
  which creates a modular group representation
  based on the data $(\Pi,p)$.
    
  \item  In \S 5, I show that
    the map $\rho$ is injective
  on $\mathcal P$ and two-to-one on ${\mathcal A\/}-{\mathcal P\/}$.
  More precisely, I consider the two components
  ${\mathcal A\/}_+$ and ${\mathcal A\/}_-$ of
  ${\mathcal A\/}-{\mathcal P\/}$. I show that
  $\rho$ is injective on each of
  ${\mathcal A\/}_+$  and ${\mathcal A\/}_-$, and I show that
$$\rho({\mathcal P\/} \cup {\mathcal A\/}_+)=\rho({\mathcal P\/} \cup {\mathcal A\/}_-)=\rho({\mathcal A\/}).$$
So, in short, $\rho$ is a kind of folding map,
just as in the real hyperbolic case.
At the end of \S 5, I do a  calculation
related to the gradient of the trace function on $\mathcal R$.
At this point, I arbitrarily choose ${\mathcal A\/}_+$ over
${\mathcal A\/}_-$.  The calculations in this chapter are
done in Mathematica [{\bf W\/}].

\item In \S 6, I work out the topology of ${\mathcal P\/} \cup {\mathcal A\/}_+$
  and show  that  $\rho: {\mathcal P\/} \cup {\mathcal A\/}_+ \to \mathcal R$
  is continuous, injective, and proper.
  Combining this with the trace calculation at the end of \S 5, I
  conclude that $\rho(\mathcal A)=\cal B$, the component
  from Theorem \ref{one}.  Hence every group in
  ${\cal A\/}-\cal P$ is Anosov and every group
  in $\cal B$ is a prism group.
  
\item In \S 7, I show that the patterns of geodesics in
  Theorem \ref{extra} are always embedded, and so
  are the patterns of flats that contain them.
\end{itemize}

One can download the Mathematica code for this paper
from \newline
{\bf http://www.math.brown.edu/$\sim$res/PappusCalcs.TAR\/}
\newline

I  thank
Martin Bridgeman,
Bill Goldman,  Tom Goodwillie,
Sean Lawton,
Joaquin Lejtreger,
Joaquin Lema, Dan Margalit,  Max
Riestenberg, Dennis
Sullivan, and Anna Wienhard for
interesting and helpful conversations.
I especially thank Martin and the
two Joaquins for recently rekindling
my interest in this subject.

      \section{The Classic Case}
\label{classic}

\subsection{The Hyperbolic Plane}

We work with the upper half-plane model $\H^2$
of the hyperbolic plane.  In this model, the
geodesics are either arcs of semicircles with
endpoints on $\R$ or else vertical rays.
The group ${\rm Isom\/}(\H^2)$ is generated
by real linear fractional transformations
and the map $z \to -\overline z$, which
is reflection in the $Y$-axis.

A group $\Lambda \subset {\rm Isom\/}(\H^2)$ acts
{\it discretely\/} if for any compact $K \subset \H^2$ there
are only finitely many $g \in \Lambda$ such that
$g(K) \cap K \not = \emptyset$.  This kind action is
also called {\it properly discontinuous\/}.
The {\it limit set\/} of $\Lambda$ is the accumulation
set on $\R \cup \infty$ of any orbit. The definition does
not depend on the orbit chosen.

\subsection{The Farey Triangulation}

The geodesics of
the Farey triangulation limit on rational
points in the ideal boundary
$\R \cup \infty$; two
rationals $a/b$ and $c/d$ are endpoints
of a geodesic in the triangulation if and only
if $|ad-bc|=1$.  See Figure 2.1

 \begin{center}
\resizebox{!}{1.7in}{\includegraphics{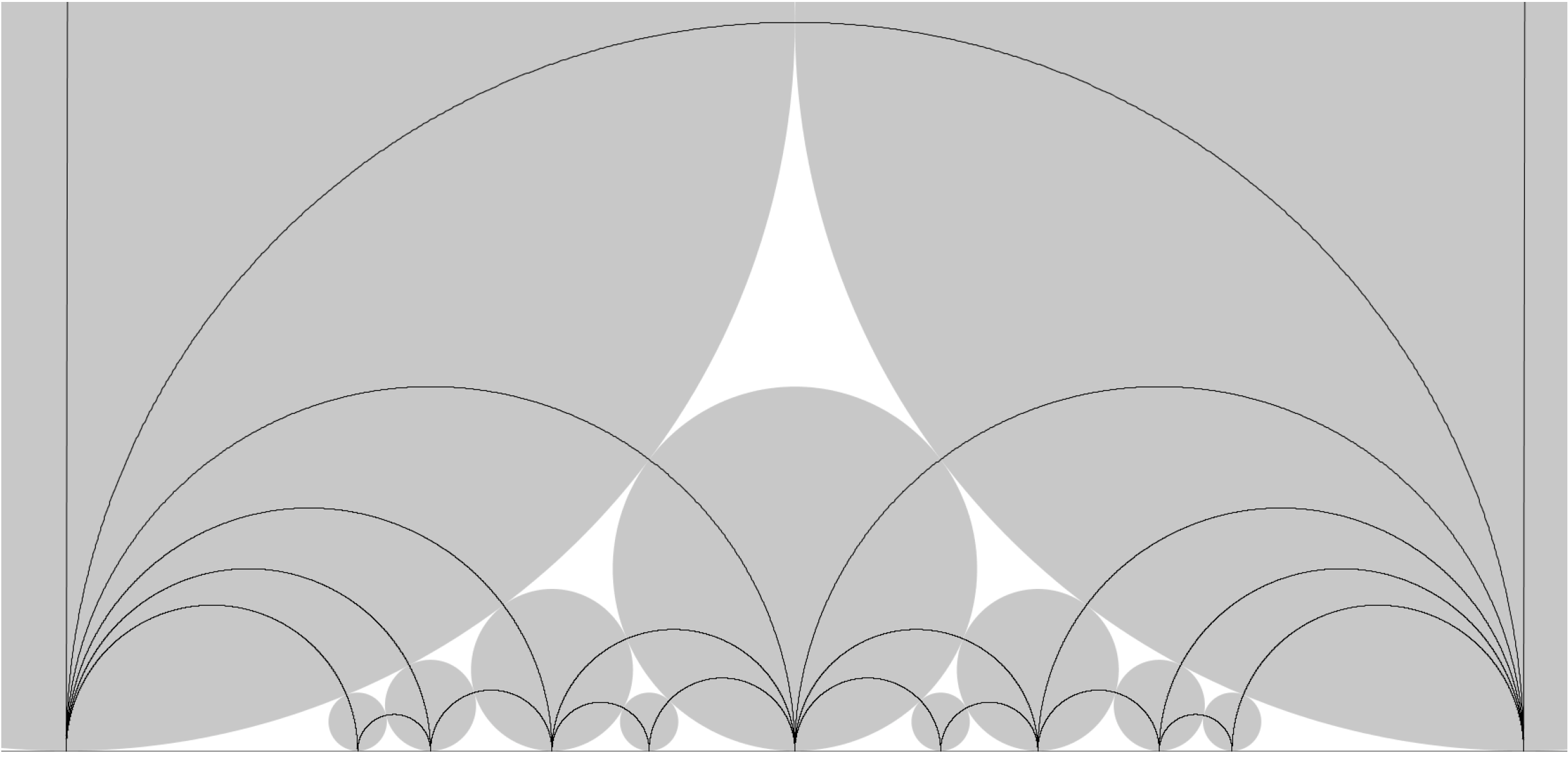}}
\newline
{\bf Figure 2.1:\/} Part of the Farey triangulation and dual horodisk
packing
\end{center}

The tangency points of the horodisks in the packing are
distinguished points on the geodesics of the Farey pattern.
We call these the {\it inflection points\/} of the geodesic.
A more robust definition of the inflection points goes like
this: Any ideal triangle in $\H^2$ has an order $6$
symmetry group.  The elements of order $2$ are
reflections about geodesics which connect an inflection
point to the opposite cusp.  This definition is nicer
because it only depends on the individual ideal
triangle.  When the triangles are arranged as in the
Farey triangulation, the robust definition coincides
with the special definition given in terms of the
horodisks.

\subsection{The Modular Group and its Shears}
\label{hypshear}

The modular group $PSL_2(\Z)$ is generated by the order
$3$ isometric rotations about the centers of the ideal
triangles in the triangulation and the order $2$ reflections
about the inflection points.  Algebraically, the modular
group is the free product $\Z/2 * \Z/3$.

The robust definition of the inflection points gives us another
way to define the modular group.  Let $\tau$ be some fixed
ideal triangle.  Then the modular group is the isometry
group generated by the order $3$ counterclockwise
rotation $\sigma_3$ of
$\tau$ and by the order $2$ rotation $\sigma_2$ about one of the
inflection points of $\tau$.  If we choose $\tau$ to
be (say) the triangle with vertices $0,1,\infty$ then
we recover the modular group exactly.  If we start with
a different choice of $\tau$ we get a group that is
conjugate to the modular group.

Now we consider shearing of the modular group.
Let $\tau$ be an ideal triangle, as above.
Let $\gamma$ be one of the geodesics
comprising $\partial \tau$.
We choose one of the  points $p\in \gamma$ which is
$s$ units from the inflection point on $\gamma$.
We then let $\Gamma_s$ denote
the group generated by $\sigma_3$ above and
by the order $2$
rotation about $p$.  When $s=0$ we recover
the modular group. When $s>0$ we get a
{\it shearing\/} of the modular group.
The other choice of $p \in \gamma$ that
is equidistant from the inflection point
gives a conjugate group.   So, the distance
$s$ here is all that really matters.

When $s \not =0$, the group $\Gamma_s$ preserves a tiling
of a closed subset $\Delta_s \subset \H^2$ by ideal
triangles.
We get this
tiling starting with $\tau$ and using the isometries to
successively lay down isometric copies of $\tau$.  Two
adjacent ideal triangles $\tau_1$ and $\tau_2$ are related
in the following way.  Let $\gamma$ be the geodesic common to
$\tau_1$ and $\tau_2$.  Then the distance between the
inflection point of $\tau_1$ on $\gamma$ and the
inflection point of $\tau_2$ on $\gamma$ is $2s$.
Thus, we can also think of getting the pair
$(\tau_1,\tau_2)$ by starting with two
adjacent triangles in the Farey triangle, sliding
one of them $2s$ units relative to the other, then
moving the union into some new position by an
isometry.

The limit set $\Lambda_s$ of $\Gamma_s$ is a
Cantor set when $s \not =0$.  The region
$\Delta_s$ is the convex hull of $\Lambda_s$.
The group $\Gamma_s$ is a classic example of an
Anosov group.  As $s \to 0$, the region
$\Delta_s$ converges to all of $\H^2$ (assuming
we keep the initial triangle $\tau$ the same for all $s$)
and the limit set $\Lambda_s$ converges to
$\R \cup \infty$  in 
the Hausdorff topology.

The description above is quite well known.
See e.g. [{\bf T\/}] or [{\bf P\/}].

\subsection{The Representation Variety}

To get a representation of the abstract
modular group $\Z/2 * \Z/3$ all we need
to do is choose an order $2$ element and
an order $3$ element.  We will insist that
our representation is in $PSL_2(\R)$, the
index $2$ subgroup of linear fractional
transformations.  We consider two
representations equivalent if they are
conjugate.   Once we do this, the only data
that is important is the distance $d$ between the
fixed points. This means that the representation
variety is just the ray $R=[0,\infty)$.
The point $0 \in R$
corresponds to the case when both
fixed points coincide.

Let $d_0$ denote the distance between the
center of an ideal triangle and its inflecton point.
Our distance parameter $d$ in this section is
related to the shear parameter $s$ from the
previous section by a formula one can work
out with hyperbolic trigonometry.  The case
$d=d_0$ corresponds to $s=0$ and the case
$d>0$ corresponds to $s \not =0$.
As is well known, the point
$d \in R$ gives rise to a discrete and
faithful (i.e. injective) representation if and only if
$d \geq d_0$.

Here is another way to describe the trichotomy.
Let
\begin{equation}
  g=\sigma_2 \circ \sigma_3.
\end{equation}
Then $g^2$ is elliptic, parabolic, or loxodromic
according as $d<d_0$, $d=d_0$, or $d>d_0$.
Since $g^2$ is given by a linear fractional
transformation based on an element of $SL_2(\R)$,
this criterion can be expressed
in terms of ${\rm Trace\/}(g^2)$ being
either less than, equal to, or greater than $2$.

Let us reconcile this with our picture of the shears.
The shears of the modular group give what
looks like a $1$-parameter family of representations
that is diffeomorphic to $\R$.  After all, we are
free to slide the point anywhere along the
geodesic $\gamma$.   However, this copy of
$\R$ maps into $R$ with a fold at $d_0$:
The image is the ray $[d_0,\infty)$. If we
shear the same amount in opposite directions,
we get conjugate groups.  In particular, every
group aside from the modular group has
$2$ distinct descriptions in terms of
shearing.  This kind of folding
picture will generalize to the case of
$X=SL_3(\R)/SO(3)$.

      \section{Geometric Preliminaries}

\subsection{The Symmetric Space}

Here I give an abbreviated account of the
corresponding material in [{\bf S1\/}].
This material is, of course, well known.
\newline
\newline
{\bf Basic Definition:\/}
The symmetric space
$X=SL_3(\R)/SO(3)$ can be interpreted as the space of unit
volume ellipsoids centered at the origin of
$\R^3$.   There is a natural origin of
$X$, the point which names
the round ball.  The group
$SL_3(\R)$ acts on $X$ in the obvious way.
If $E$ is an ellipsoid and $T \in SL_3(\R)$ then
$T(E)$ is just the image ellipsoid. Here I am
somewhat blurring the distinction between
points in $X$ and the ellipsoids they name.
The group $SL_3(\R)$ acts transitively on
$X$ and the stabilizer of the origin is
$SO(3)$.  So, the orbit map gives an
isomorphism between the coset description of
$X$ and the ellipsoid description.

One can also interpret $X$ as the space of
unit determinant positive definite symmetric
matrices.  Each matrix like this defines an
inner product on $\R^3$ and this unit ball
of this inner product is a unit volume ellipsoid
centered at the origin.  This is how the
correspondence between two descriptions works.
   \newline
\newline
{\bf The Metric:\/}
The space $X$ has a canonical
$SL_3(\R)$ invariant metric which
is induced by a Riemannian metric of
non-positive sectional curvature.
The distance between $E_0$ and
the {\it standard ellipsoid\/} $E(a,b,c)$ given by
\begin{equation}
\frac{x^2}{a^2}+\frac{y^2}{b^2}+\frac{z^2}{c^2}=1, \hskip 30 pt
a,b,c>0, \hskip 20 pt abc=1.
 \end{equation}
is
\begin{equation}
  \|(\log(a),\log(b),\log(c))\|=
  \sqrt{\log^2(a)+\log^2(b)+\log^2(c)}.
\end{equation}
The rest of the metric can be deduced from
symmetry.
\newline
\newline
{\bf Isometries:\/}
As already mentioned, $SL_3(\R)$ acts isometrically on $X$.
There is also an order $2$ isometry $\Delta$ of $X$ which fixes
the origin and reverses all the geodesics through the origin.
In terms of the matrix interpretation of $X$, this isometry
is given by $S \to S^{-1}$, where $S$ is a positive definite
symmetric matrix.   We call this map the
{\it standard polarity\/}, for reasons discussed in \S \ref{connection}.
The standard polarity maps the ellipsoid
$E(a,b,c)$ to $E(1/a,1/b,1/c)$.  

The group ${\rm Isom\/}(X)$ is generated by
$\Delta$ and $SL_3(\R)$.  We call the elements
of $SL_3(\R)$ {\it projective transformations\/} and
the remaining elements {\it dualities\/}.  In particular, any
point of $X$ is fixed by an order $2$ isometry
(a conjugate of $\Delta$) which reverses all
the geodesics through that point. We will call such an
isometry an  {\it elliptic polarity\/}.   These names are
usually reserved for maps on the projective plane and
its dual, but here we are talking about the action on $X$.
We explain the connection in \S \ref{connection}.
\newline
\newline
{\bf Flats:\/}
The {\it standard flat\/} $F_0$ is the union of all the points
representing standard ellipsoids.  The rank $2$
abelian group of diagonal matrices acts transitively
on the standard flat.  Thus, $F_0$ is
isometric to a Euclidean plane.
In particular, the straight lines in $F_0$ are geodesics in $X$.
Every other flat in $X$ is isometric to $F_0$.  In particular,
the structure of $F_0$ determines the structure of all the flats.

There are $3$ {\it singular geodesics\/} through the origin
in $F_0$: These correspond to the standard ellipsoids
where the set $\{a,b,c\}$ has cardinality at most $2$.
That is, either $a=b$ or $a=c$ or $b=c$.   In general,
the {\it singular geodesics\/} in $F_0$ are the ones
parallel to the singular geodesics through the origin.
A geodesic in $F_0$ is contained in more than one
flat if and only if it is a singular geodesic.
All other geodesics in $F_0$ lie only in $F_0$.

There are $3$ {\it medial geodesics\/} though
the origin. These correspond to triples $(a,b,c)$ where
either $a=1$ or $b=1$ or $c=1$.  More generally,
a {\it medial geodesic\/} in $F_0$ is one parallel to
a medial geodesic through the origin.
Each medial geodesic lies in a unique flat.
In terms of the cyclic order on the geodesics
through the origin in $F_0$, the singular
geodesics alternate with the medial geodesics,
and the angle between adjacent singular
and medial geodesic is $\pi/6$.   
A {\it medial foliation\/}
of a flat is a foliation by parallel medial geodesics.
Thus, every flat has $3$ medial foliations.  We call
a flat with a distinguished medial foliation a
{\it marked flat\/}.
\newline
\newline
{\bf Visual Boundary:\/}
The visual boundary of $X$ is defined to be the union
of geodesic rays through the origin.
We denote this as $\partial X$.  The action of
isometries on $X$ extends to give a homeomorphism
of $\partial X$ in the following way. If
$\rho$ is a geodesic ray through the origin and $I$ is
an isometry then the image of $\rho$ under $I$ is some
other geodesic ray, not necessarily contained on a geodesic
through the origin.  There is a unique geodesic ray through
origin $\rho'$ such that the distance between corresponding
points of $I(\rho)$ and $\rho'$ remains uniformly bounded.
The action of $I$ on $\partial X$ maps $\rho$ to $\rho'$.

\subsection{Connection to Projective Geometry}
\label{connection}

{\bf Projective Objects:\/}
The projective plane $\P$ is the set of
$1$-dimensional subspaces of $\R^3$. 
The dual plane $\P^*$ is the set of
$2$-dimensional subspaces of $\R^3$.
The {\it flag variety\/} is the set of
pairs $(p,\ell)$ where $p$ is a
$1$-dimensional subspace of $\R^3$
and $\ell$ is a $2$-dimensional subspace
of $\R^3$, and $p \subset \ell$.
These objects are called {\it flags\/}.
Equivalently, a flag is a pair
$(p,\ell)$ where $p$ is a point of
$\P$ and $\ell$ is a line of $\P$,
and $p \in L$.   Each point in
$\P$ corresponds to a line in
$\P$, namely the set of $1$-dimensional
subspaces contained in a given
$2$-dimensional subspace.
\newline
\newline
{\bf Limits of Singular and Medial Geodesics:\/}
The singular geodesics accumulate at one end
to points of $\P$ and at the other end to
points of $\P^*$.   This is easily seen for
the standard flat.  For one of the singular
geodesics through the origin, the corresponding
standard ellipsoids are $E(a,a,1/a^2)$.  As
$a \to 0$ these become long and thin and
pick out a $1$-dimensional subspace in $\R^3$.
As $a \to \infty$,  these ellipsoids flatten
out like a pancake and define a $2$ dimensional
subspace of $\R^3$.

The medial geodesics accumulate at both
end at points of the flag variety.  The standard
example is the medial geodesic consisting
of $E(a,1,1/a)$.  As $a \to 0$ the $1$-dimensional
subspace is the $Z$-axis and the $2$-dimensional
subspace is the $YZ$-plane. As $a \to \infty$
the $1$-dimensional subspace is the $X$-axis
and the $2$-dimensional subspace is the
$XY$-plane.  The intuition here is that in either
direction these ellipsoids look like popsicle
sticks.   The longest direction picks out the
one dimensional subspace and the two longest
directions pick out the two dimensional subspace.
\newline
\newline
{\bf Marked Flats and Pairs of Flags:\/}
A triple of points in $\P$ is in {\it general position\/}
if they are not contained in the same line.
Likewise, a triple of lines in $\P$ is in
general position if they do not have a single point
in common.
Two flags $(p_1,\ell_1)$ and $(p_2,\ell_2)$ are
{\it transverse\/}  if $p_1 \not \in \ell_2$ and
$p_2 \not \in \ell_1$.

A marked flat defines
a pair of transverse flags, namely
the (common) limits of the medial geodesics in the
foliation.
Conversely a pair of transverse flags
determines a unique marked flat. The space
of marked flats is $6$-dimensional.
A {\it projective triangle\/}, namely
a triple of general position points,
determines a unique flat, and
{\it vice versa\/}.  A pair of
transverse flags determines a
unique projective triangle, and
a projective triangle determines
three pairs of transverse flags,
corresponding to the three markings of the flat.
\newline
\newline
{\bf Our Terminology Explained:\/}
Consider an element of ${\rm Isom\/}(X)$ which
we are calling a projective transformation.
Such an element acts on $\partial X$ in such a
way as to preserve $\P$.  The action the map
induces on $\P$ is precisely that of a classical
projective transformation: It is a homeomorphism
of $\P$ which maps lines to lines.

Consider an element
of ${\rm Isom\/}(X)$ which we are calling
a duality.  Such an element acts on
$\partial X$ in such a way as to
swap $\P$ and $\P^*$.  The induced
map we get on $\P \cup \P^*$ is a
duality in the classic sense: It swaps points
and lines in such a way as to carry
collinear lines to coincident lines.

In particular, the isometries of $X$ we are calling
{\it elliptic polarities\/} induce the
maps on $\P \cup \P^*$ which are
commonly called elliptic polarities.
Such maps are given by
$\ell \leftrightarrow \ell^*$ where the $1$-dimensional
subspace
in $\R^3$ representing $\ell$ is
mapped to the $2$-dimensional
subspace representing $\ell^*$ which
is orthogonal to it with respect to
some inner product.
When this inner product is the
usual dot product, we get the
standard polarity.

\subsection{Matrix Actions}

Here we explain how we compute
the action of projective transformations
and dualities using matrices.
\newline
\newline
{\bf Representing Points and Lines:\/}
We represent points in $\P$ as $3$-vectors.
When $c \not =0$, the
vector $(a,b,c)$ represents the point
$(a/c,b/c)$ in the {\it affine patch\/}.
The affine patch is essentially a copy of
$\R^2$ sitting inside $\P$.
We also represent lines as vectors.
The vector $(a,b,c)$ represents the
line given by the subspace $ax+by+cz=0$.
If we have two vectors $v_1=(a_1,b_1,1)$ and
$v_2=(a_2,b_2,1)$ then the vector
$(1+t) v_1 - t v_2$ represents a point on
the line $\overline{v_1v_2}$.
\newline
\newline
{\bf Action on Points:\/}
We will work with matrices in $GL_3(\R)$.
For our purposes we do not need to fuss about
whether our matrix has determinant $1$.  We
can always scale the matrix to have this property.
(Here we are strongly using the property that
we have $n \times n$ matrices and $n=3$ is odd;
in the even case the situation is a bit different.)
The matrix $S$ acts on a vector $\widehat p$ representing
a point $p$ by linear transformation:  The new vector
$S(\widehat p)$ represents $S(p)$.
\newline
\newline
{\bf Action on Lines:\/}
The matrix $S$ acts on our line representations
in the following way:  We let
$(S^{-1})^t$ act on the vector representation $\widehat L$
of a line $L$.
and the new vector represents the line $S(L)$.
Here we are taking the inverse-transpose.
A few calculations will convince the reader that this
is indeed the right thing to do.   This works because
$$S(\widehat p) \cdot (S^{-1})^t(\widehat L )=
S^{-1}S(\widehat p) \cdot \widehat L=\widehat p \cdot \widehat L.$$
This way of defining the action is correct because it preserves
incidences between points and lines.
\newline
\newline
{\bf Dualities:\/}
The standard polarity $\Delta$ just acts as the
identity matrix, both on points and lines.  Thus
$\Delta(a,b,c)=(a,b,c)$. All that changes is the
interpretation of the meaning of the vector $(a,b,c)$.
It is most convenient to represent dualities as
compositions $\Delta \circ S$.   We have the equations
\begin{equation}
S \circ \Delta = \Delta \circ (S^{-1})^t , \hskip 30 pt
S \circ \Delta \circ S^{-1}= \Delta \circ (S^{-1})^t S^{-1}.
\end{equation}
The first of these equations implies the second one.
\newline
\newline
{\bf The Representation Space:\/}
The modular group $G=\Z/2*\Z/3$ is generated by
$\sigma_2$ and $\sigma_3$, elements of order $2$ and $3$.
We consider representations of $G$ into ${\rm Isom\/}(X)$ such that
\begin{itemize}
\item $\rho(\sigma_2)$ is an elliptic polarity.
\item  $\rho(\sigma_3)$ is the matrix given by
\begin{equation}
  \label{rotation}
  \left[\matrix{\cos(\theta)&\sin(\theta)&0 \cr
      -\sin(\theta)&\cos(\theta)&0 \cr 0&0&1}\right], \hskip 30 pt
  \theta=2 \pi/3.
\end{equation}
This matrix acts on
$\P$ so as to preserve the affine patch and
act there as an order $3$ rotation about the origin.
\item The fixed point of $\rho(\sigma_2)$ does not
  lie in the fixed point set of $\rho(\sigma_3)$.
\end{itemize}
The fixed point set
of $\rho(\sigma_3)$ is the singular geodesic $\gamma$.
consisting of standard ellipsoids $E(a,a,a^{-2})$.
We call these representations
{\it normalized\/}.
We consider two representations to be
the same if
they are conjugate in ${\rm Isom\/}(X)$.

Let $\cal R$ denote the space of all normalized
representations, modulo conjugacy.
The distance between two elements
$[\rho_1], [\rho_2] \in \cal R$
to be the minimal $D$ such that there are
two normalized representatives
$\rho_1$ and $\rho_2$ such that
the fixed point sets of $\rho_1(\sigma_1)$
and $\rho_2(\sigma_2)$ are $D$ apart in $X$.
In [{\bf S2\/}] we proved the following result.

\begin{theorem}
  \label{trace}
  $\cal R$ is homeomorphic to $\R^3-\{(0,0,0)\}$
  and is a smooth manifold away
    from the two curves, one corresponding to
    line-preserving representations and one corresponding
    to conic-preserving representations.
   The trace  of any word is a smooth function on the smooth
  points of $\cal R$.
\end{theorem}

      \section{The Prism Representations}

\subsection{Basic Defintions}

We will study triples of flags $\{(p_i,\ell_i)\}$.
We insist that $p_i \not \in \ell_j$ when $i \not = j$. 
Let $P_1,P_2,P_3$ and $L_1,L_2,L_3$ respectively be
vectors representing $p_1,p_2,p_3$ and $\ell_1,\ell_2,\ell_3$. 
These vector representatives are unique up to scaling.
Our condition is that $P_i \cdot L_j \not =0$ when $i \not = j$.
The {\it triple product\/} of our flag triple is
\begin{equation}
  \label{tripp}
 \chi= \frac{(P_1 \cdot L_2)(P_2 \cdot L_3)(P_3 \cdot L_1)}{(P_2 \cdot L_1)(P_3 \cdot L_2)(P_1\cdot L_3)}.
\end{equation}
This is a very well known invariant.
Compare [{\bf FG\/}].

We call the triple of flags {\it negative\/} if
$\chi < 0$.
We can always normalize by a projective transformation
so that our triple of flags is rotationally symmetric.
Figure 4.1 shows the $4$ possibilities. 
In the first two cases the points are not visible because
they lie in the line at infinity.

\begin{center}
\resizebox{!}{1.3in}{\includegraphics{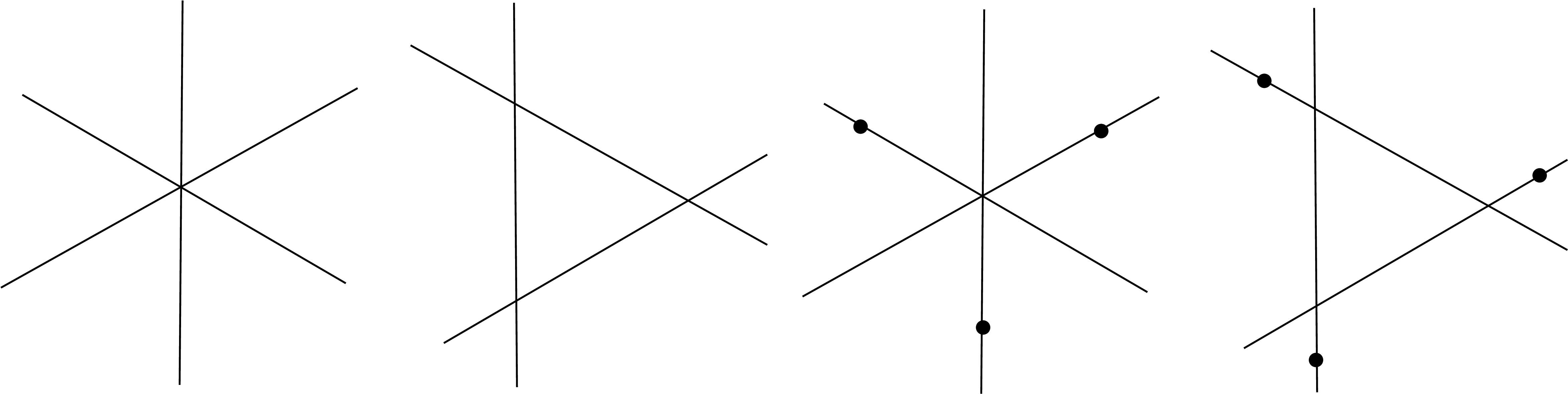}}
\newline
{\bf Figure 4.1:\/}  Negative triples with
$3$-fold Euclidean symmetry.
\end{center}

The figure on the right is meant
to represent an infinite family of examples.
In the figure on the right,
we are allowed to move the points in and out in a symmetric
way as long as we do not move them into the central triangle.
It is worth remarking, finally, that the
middle two cases are dual to each other.

We call the triple of flags {\it generic\/} if it is projectively
equivalent to a member of the family represented by the
figure on the right in Figure 4.1.

If we permute the order of our flags then $\xi$ is either
preserved or replaced by $1/\xi$.
For non-generic negative triples we have $\xi=-1$ and the invariant
cannot tell apart the various cases.  For the generic case, two
negative triples are projectively equivalent if and only if
they have the same triple invariant.   Referring to Figure 4.1,
the invariant $\chi$, when not equal to $-1$,
measures how the triple of points is placed with respect
to the inner equilateral triangle.
\newline
\newline
\noindent
{\bf Definition:\/}
A {\it prism\/} is the triple of marked flats
corresponding to a negative triple of flags. 
We define the {\it triple invariant\/} of the prism $\Pi$ to be
\begin{equation}
 \xi(\Pi)= |\log(-\chi)| \in [0,\infty),
\end{equation}
where $\chi$ is the triple invariant of a triple
of flags defining $\Pi$.  We take absolute
values so as to get an invariant that is independent
of permutations of the flags and also self-dual.
By construction, the
prism $\Pi$ is generic if and only if
$\xi(\Pi)>0$.

\begin{lemma}
  The generic prisms $\Pi_1$ and $\Pi_2$
  are isometric iff  $\xi(\Pi_1)=\xi(\Pi_2)$.
\end{lemma}

\startproof
Suppose $\Pi_1, \Pi_2$ are generic
prisms and $I: P_1 \to P_2$ is an
isometry.  Then $I$ maps the one triple of
flags to the other and either preserves the
triple product or inverts it (depending on
whether $I$ comes from a projective
transformation or a duality.)
Conversely, any two triples with the
same or reciprocal triple invariants are
equivalent under some isometry of $X$.
\endproof

\subsection{Inflection Points and Lines}

In this section we pick out some special geometric
features of prisms, which we call inflection points
and inflection lines. The inflection points only
exist for the generic prism and the inflection
lines exist in all cases.

\begin{lemma}
  \label{order2}
  The symmetry group of a generic prism is
  isomorphic to $S_3$, the permutation group of
  order $6$.  The even permutations
  are induced by projective transformations
  and the odd permutations are induced by polarities.
  \end{lemma}

  \startproof
  This is proved in [{\bf S1\/}].  Here is a sketch.
Suppose we apply the standard polarity to our flag
triple.  We then get the same vector representatives
except that their roles have changed.   Therefore,
the triple product of the dual triple is the reciprocal
of the original triple product. If we then apply an odd permutation
to the flags we get back to the original invariant.
This operation implies the existence of an
order $2$ symmetry of the flag, induced by
a polarity, which does an odd permutation to
the flats comprising the prism.

The $2$-fold symmetry just explained combines with the
$3$-fold symmetry to give us a symmetry
group $H=S_3$ of order $6$.   Suppose $\psi$ is
some other symmetry. Composing with
some element of $H$ we can consider
the case when $\psi$ preserves at least
one flag and also is a projective transformation.
But then $\psi$ has to induce the identity
permutation on the flags because of the
triple invariant.  But then $\psi$ is a projective
transformation which fixes $6$ general position
points.  Hence $\psi$ is the identity.
This shows that $H$ is the full group of symmetries.
\endproof

\begin{lemma}
  Let $\Pi$ be a generic prism.  Each order $2$ isometry
  of $\Pi$ fixes a unique point in the flat of $\Pi$ that
  it stabilizes.
\end{lemma}

\startproof
This is proved in [{\bf S1\/}].  Here is the proof again.
Let $\delta$ be such an isometry and let $F$ be the
flat such that $\delta(F)=F$.  The duality
$\delta$ swaps the two flags defining $F$ and
hence reverses the directions of the medial geodesics
foliating $F$ and asymptotic to these flags.  Also,
being a polarity, $\delta$ reverses the directions of
all singular geodesics in $F$.  In particular
$\delta$ reverses an orthogonal pair of
directions. This forces $\delta$ to reverse
every direction.   If we identify $F$ with
$\R^2$ then $\delta$ is acting as an isometry
whose linear part is an order $2$ rotation.
Such maps have unique fixed points in $\R^2$.
\endproof

\noindent
{\bf Definition:\/}
On a generic prism $\Pi$, the
{\it inflection points\/} are
the fixed points of the order
$2$ isometries of $\Pi$.
There are $3$ such inflection
points, one per flat.   They
are permuted by the order
$3$ isometries of $\Pi$.
The {\it inflection lines\/} are
the singular geodesics which
contain the inflection points and
which are perpendicular to the
geodesics in the medial foliations.
\newline

For a prism based on either of the
two middle pictures in Figure 4.1, the
inflection points do not exist.
Geometrically, what is happening as
we approach one of these prisms
through a family of generic prisms is that
the inflection points move
off to $\infty$.   The inflection
lines still exist however, as
we now explain.

In [{\bf S1\/}]
we show that every Pappus modular group
is an isometry group of an embedded
pattern of flats.   In the generic case
we show that each fixed point of an
order $2$ element of the group is
contained in the relevant inflection line.
(See also \ref{inflectionproof}.)
If we exclude the totally symmetric
Pappus modular groups, the remaining
$1$-parameter family of non-generic
groups can be normalized so that they
all involve the same flat.  Taking
a limit of the generic result, we can say
that all the order $2$ fixed points of
all these groups in $F$
lie on the same singular geodesic which
is perpendicular to the medial foliation of $F$.
This singular geodesic is the inflection line in $F$.

We have not yet discussed the totally
symmetric case, the prism based on the
lefthand picture in Figure 4.1.
The associated prism has an infinite symmetry group.
Referring to Figure 1, the projective
transformation which extends the map
$x \to rx$, for any $r \not =0$, induces
an isometry that preserves the prism.
These isometries act nontrivially on the flats.
In this case every associated triangle is isometric
to a hyperbolic Farey triangle.  The inflection
lines are comprised of the symmetry points on
each ideal hyperbolic triangle.

\subsection{Triangle Foliations}

Let $\Pi$ be a prism.
The order $3$ isometries of $\Pi$ preserve
the medial geodesic foliations.
Thus $\Pi$ is foliated by
{\it triangles\/}, triples of medial
geodesics invariant under the order
$3$ isometries of $\Pi$.
In the generic case, exactly one
triangle of $\Pi$ contains all
$3$ inflection points. In all
cases, the triangles of $\Pi$
are perpendicular to the inflection lines.
In particular, this fact gives us a way to
pick out canonical {\it inflection points\/}
on each triangle, namely where the
triangle intersects the inflection line.
In the totally symmetric case,
all the triangles are isometric to hyperbolic ideal
triangles and hence isometric to each other.
For the other prisms the situation is
very different.

\begin{lemma}
  \label{isom}
  Let $\Pi_1$ and $\Pi_2$ be generic prisms
  and let $\tau_k$ be a triangle of $\Pi_k$
  for $k=1,2$.   Then
  $\tau_1$ and $\tau_2$ are
  isometric to each other only if
    $\Pi_1$ and $\Pi_2$ are isomorphic prisms.
    If $\Pi_1=\Pi_2=\Pi$ then $\tau_1$ and $\tau_2$ are
isometric if and only if they are permuted by
the symmetry group of $\Pi$.
\end{lemma}

\startproof
An isometry taking $\tau_1$ to $\tau_2$ would
have to map the flats of $\Pi_1$ to the flats
of $\Pi_2$.  This proves the first statement.
For the second statement, note that a
projective symmetry of $\Pi$ taking
$\gamma_1$ to $\gamma_2$ must be
in the symmetry group of $\Pi$.
\endproof

\subsection{The Axis}
\label{axis}

In this section we prove a properness result
about prisms that will come in handy when we analyze
components of the representation variety.  We
first need to define what we mean the axis
of a prism.

\begin{lemma}
  The fixed point set of the
  order $3$ symmetries of a prism is a
  singular geodesic.
\end{lemma}

\startproof
This is a general fact about order $3$ isometries.
We can normalize so that the order $3$
symmetries are as in Equation \ref{rotation}.
The fixed point set is the singular
geodesic corresponding to the standard
ellipses of the
form $E(a,a,a^{-2})$.
\endproof

We call the fixed point set of the order $3$ symmetry the
{\it axis\/} of the prism.

Our next result compares two geometric
properties of prisms.  The second of
these quantities is related to the topology
of the representation space $\cal R$.
The two quantities compared here are like
our two distance parameters $s$ and $d$
discussed in \S \ref{classic}.
Given a prism $\Pi$ and a point $p\in \Pi$.
we define $s(\Pi,p)$ be the distance
from $p$ to the inflection line in the flat
of $\Pi$ that contains $p$.  At the same
time, let $d(\Pi,p)$ be the distance
from $p$ to the axis of $\Pi$.

\begin{lemma}[Properness]
  \label{proper}
  Let $\{\Pi_n,p_n\})$ be any sequence of
  prisms.  If  $s(\Pi_n,p_n) \to \infty$
  then also $d(\Pi_n,p_n) \to \infty$.
\end{lemma}

\startproof
We will suppose this false and derive a contradiction.
That is, we suppose that
$s(\Pi_n,p_n) \to \infty$ but
$d(\Pi_n,p_n)$ stays bounded.
We can normalize by isometries so that
the point on the axis of $\Pi_n$ closest
to $p_n$ is the origin of $X$.  This
means that the distance from $p_n$ to
the origin is uniformly bounded.
But then the flat $F_n$ of $\Pi_n$ containing
$p_n$ intersects a uniformly bounded
region of $X$.   Since the order $3$ isometries
of $\Pi_n$ fix the origin, we see that all
flats of $\Pi_n$ intersect a uniformly
bounded region of $X$.

But then we can take a limit and get
a prism $\Pi=\lim \Pi_n$.
Since the inflection lines of $\Pi$
exist and are unique, we see that the
inflection lines of $\Pi_n$ remain
within a uniformly bounded region
of $X$.   But then we have a uniformly
bounded distance from $p_n$ to the relevant
inflection line, a contradiction.
\endproof

\subsection{Modular Group Representations}
\label{PRISM}

We say that a {\it prism pair\/} is a
pair $(\Pi,p)$ where $\Pi$ is a prism
and $p \in \Pi$.  We mean here that
$p$ is a point in one of the flats comprising $\Pi$.
We impose a cyclic
order on $\Pi$, determined by the cyclic
order on the flats.  The order $3$
symmetries of $\Pi$ respect this order
and the order $2$ symmetries do not.
One of the order $3$ symmetries cycles
the flats of $\Pi$ one click forward in the
cyclic order and the other one cycles
the flats of $\Pi$ one click backward.
We prefer the former symmetry and
we call it the {\it forward symmetry\/}.
When we normalize as in Figure 4.1, the
forward symmetry is given by the
element $\rho(\sigma_3)$ from
Equation \ref{rotation}

The prism pair $(\Pi,p)$ determines a
point in $\cal R$.  We let $\rho(\Pi,p)$
be the representation such that
$\rho(\sigma_2)$ is the elliptic
polarity fixing $p$ and $\rho(\sigma_3)$ is
the forward symmetry of $\Pi$.
We call these the
{\it prism representations\/}.

We call two prism pairs $(\Pi_1,p_1)$ and
$(\Pi_2,p_2)$ {\it equivalent\/} if there is an
isometry of $X$ which maps the
first pair to the second and respects
the imposed cyclic orders.
We let $\cal A$ denote the space of
equivalence classes of prism pairs.
We have a map ${\cal A\/} \to {\cal R\/}$.
Here ${\cal R\/} \cong \R^3-\{(0,0,0)\}$ is the
big representation space we considered
in the previous chapter.
We call a prism pair $(\Pi,p)$ {\it neutral\/}
if $p$ lies on an inflection line of $\Pi$.
We proved in [{\bf S1\/}] that every
neutral prism pair gives rise to a Pappus
representation of the modular group
and conversely that every Pappus
representation of the modular
group arises this way.  Also see
\S \ref{inflectionproof}.

We let ${\cal P\/} \subset {\cal A\/}$ denote the
set of neutral prism pairs.   Our results in the next
chapter will show that the map
$\rho: {\cal A\/} \to {\cal R\/}$ is
one-to-one on $\cal P$ and
two-to-one on ${\cal A\/}-{\cal P\/}$. This
  result generalizes
the folding phenomenon we discussed in
\S 2 in the hyperbolic setting.

      \section{The Big Calculation}

\subsection{The Main Results}

We continue the notation from the last
section of the previous chapter.
Given $\rho=\rho(\Pi,p)$ define
\begin{equation}
  \label{key}
 g_{\rho}=\rho(\sigma_2 \sigma_3).
\end{equation}
This element $g=g_{\rho}$ preserves one of the flags
$f_1$ associated to the flat $F_1$ of $\Pi$ that contains $p$.
To see this, let $f_1,f_2,f_3$ be the flags defining $\Pi$,
chosen so that $F_1$ is determined by the pair $(f_1,f_2)$.
Then $\sigma_3(f_1)=f_2$ and $\sigma_2(f_2)=f_1$.
Hence $g(f_1)=f_1$.
The square $g^2$ also preserves $f_1$.
It is easier to work 
with $g^2$ because this element is a projective transformation.
Below we prove the following result:

    \begin{theorem}
    \label{loxo}
    The element $g_{\rho}^2$ is parabolic iff
    $p$ lies on the inflection line of $\Pi$.  This happens
    iff $\rho$ is a Pappus modular group representation.
    Otherwise
    $g_{\rho}^2$ has eigenvalues $(\lambda,1/\lambda,1)$
    with $\lambda \in (-\infty,-1) \cup (-1,0)$.
    We can choose $\lambda$ so that the corresponding
    eigenvector corresponds to the flag $f_1$.
  \end{theorem}

  The final statement requires some explanation.
  To keep consistent with our notation below, we
  write $f_1=(b_1,L_2)$.   What we are saying,
  first of all, is that the eigenvector of $g^2$
  corresponding to $\lambda$ represents $b_1$.
  We are also saying that $\lambda$ is an eigenvalue
  of $(g^{-2})^t$, and the corresponding eigenvector
  represents $L_2$.

  We call the prism pair $(\Pi_1,p_1)$ {\it attracting\/}
  if $|\lambda|>1$. This property is independent of
  how we normalize $(\Pi,p)$.  This is obvious if we
  replace $(\Pi,p)$ by some pair
  $(T(\Pi),T(p))$ where $T$ is a projective transformation.
  This is far less obvious if we take $T$ to be a duality.
  The reader might worry that somehow $\lambda$ gets
  changed to $1/\lambda$.  This is not the case. One
  way to check this is just to try some experiments
  with diagonal matrices and the standard flags
  associated to them.  Another way is to observe
  that the attracting nature of $(\Pi,p)$ has a
  geometric interpretation in terms of the symmetric
  space $X$:  The isometry $g^2$ is moving points
  in $X$ towards the point in the visual boundary
  corresponding to $f$.  This is an isometry-invariant
  way to talk about the attracting nature of $(\Pi,p)$.
  If $|\lambda|<1$ we call $(\Pi,p)$ {\it repelling\/}.
  Finally, as in the previous chapter, we call $(\Pi,p)$
  {\it neutral\/} if $\lambda=-1$.
  
  The element $g^2$ also has an eigenvalue $1/\lambda$
  and there is some other flag $f_1'$ that corresponds to
  this eigenvalue.   below we also prove the following result.
  
    \begin{theorem}
      \label{loxo2}
      If $\rho$ is not a Pappus modular group representation, then the
     orbit of $f'_1$ under $\rho(\sigma_3)$ defines a
     prism $\Pi'$ such that
     $\rho(\Pi',p)=\rho(\Pi,p)$.
    Exactly one of the prism pairs is attracting and
    exactly one is repelling.
  \end{theorem}
  We will prove Theorems \ref{loxo} and \ref{loxo2}
  in this chapter.
  
  Recall that $\cal A$ is the space of isometry classes of
  prism pairs.  Let ${\cal A\/}_+$ denote the set of
  attracting prism pairs.  Our corollary below favors the
  attracting prism pairs over the repelling prism
  pairs, but we could make the same kind of
  statement about the repelling pairs.  Let
  $\rho: {\cal A\/} \to {\cal R\/}$ be the map
  which assigns each isometry class of prism
  pair its representation class in $\cal R$.
  We are slightly abusing notation here,
  because $\rho(\Pi,p)$ is also denoting the
  individual representation based on
  $(\Pi,p)$ and not its conjugacy class.

  \begin{corollary}
    \label{inj}
    The map $\rho$ is injective on
    ${\cal P\/} \cup {\cal A\/}_+$ and
    $\rho({\cal P\/} \cup {\cal A\/}_+)= \rho({\cal A\/})$.
  \end{corollary}

  \startproof
  Certainly a neutral prism pair cannot give the
  same representation as an attracting or
  repelling pair because parabolic elements
  are not conjugate to loxodromic elements.
  Hence $\rho({\cal P\/}) \cap \rho({\cal A\/}_+)=\emptyset$.

  Suppose $\rho(P_1,p_1)=\rho(P_2,p_2)$ for two prism pairs
  in ${\cal P\/} \cup {\cal A\/}_+$.   Such that
  $\rho(\Pi_1,p_1)=\rho(\Pi_1,p_2)$.
  (We can adjust by an isometry so that these
  representations are equal and not just conjugate.)
  The common element $g^2$ cannot be both
  loxodromic and parabolic. Hence both prism
  pairs are either neutral or attracting.

  Consider the neutral case first.
  The element $g^2$ has a unique fixed flag
  $f_1$, and $f_1$ must be one of the triple
  of flags defining both $\Pi_1$ and $\Pi_2$.
  But then the orbit of $f_1$ under $\rho_k(\sigma_3)$
  defines $\Pi_k$.  Since $\rho_1(\sigma_3)=\rho_1(\sigma_3)$
  we see that $\Pi_1=\Pi_2$.  Since
  $\rho_1(\sigma_2)=\rho_2(\sigma_2)$ and
  $p_k$ is the unique fixed point of
  $\rho_k(\sigma_2)$, we have
  $p_1=p_2$.

  Now consider the attracting case.   One of the
  flags $f_1$ defining $\Pi_k$ is the attracting fixed
  point of $g_k^2$ for each $k=1,2$.  Since
  these are the same element, the same flag
  is part of the triple defining both
  $\Pi_1$ and $\Pi_2$.  But then the orbit of
  this flag under the common element
  $\rho_k(\sigma_3)$ gives the triple
  defining $\Pi_k$.  Hence $\Pi_1=\Pi_2$.
  Likewise $p_1=p_2$.  This proves the
  first statement of the lemma.

  The second statement follows from
  Theorem \ref{loxo2}, which says that
  each non-neutral member of $\cal A$
  has the property that there are both
  attracting and repelling pairs which
  give the same representation.
\endproof

\subsection{Normalizing Triples of Flags}
\label{normal}

As preparation for proving Theorems \ref{loxo} and
\ref{loxo2} we discuss how to normalize triples of
flags. We consider the generic case, and then at
the end of our calculations consider the non-generic case.
   We can normalize the picture as in the right-hand
   picture in Figure 4.1.  Figure 5.1 repeats with this
   picture, and with labels.
   The flags are $f_k=(b_k,L_{k+1})$.  Here and
   everywhere else we take the indices mod $3$.

 \begin{center}
\resizebox{!}{2.5in}{\includegraphics{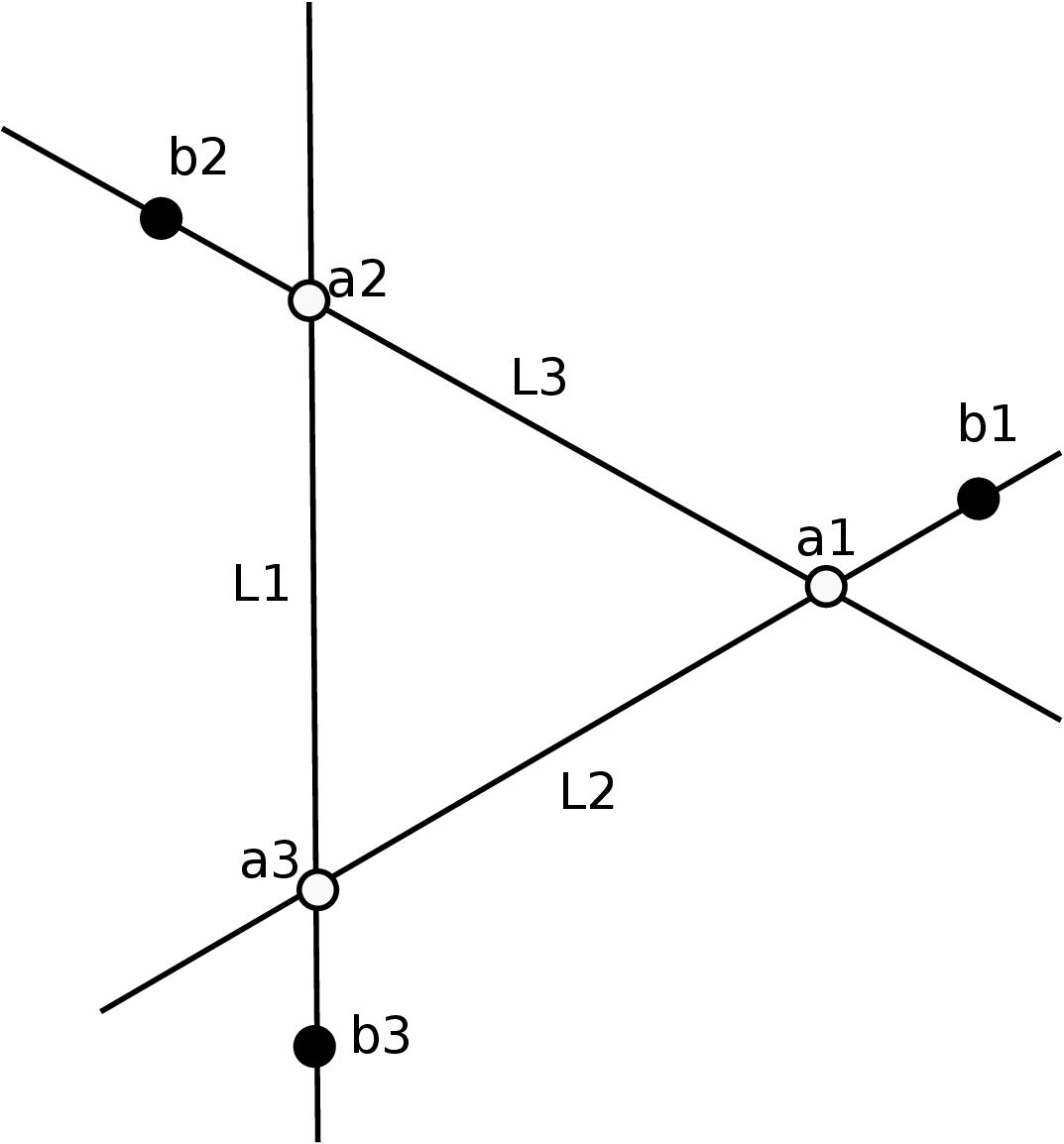}}
\newline
{\bf Figure 5.1:\/}  A normlized Flag
\end{center}

In the case shown in Figure 5.1, the point
$a_1$ is between $a_3$ and $b_1$ on the
line $L_2$.  This case corresponds to the
triple invariant of $(f_1,f_2,f_3)$ lying in
$(-1,0)$.   The other case would be when
$a_3$ lies between $b_1$ and $a_1$.
This corresponds to the triple invariant
lying in $(-\infty,-1)$.   The intermediate
case, when $b_1,b_2,b_3$ all lie on the line
at infinity, corresponds to the triple
invariant being equal to $-1$.

We can apply the standard duality $\Delta$ to the
picture.  The new flags $\Delta(f_1),\Delta(f_2),\Delta(f_3)$
have two properties we remark on:
\begin{enumerate}
\item The order $3$ counterclockwise rotation about the origin
has the action $\Delta(f_k) \to \Delta(f_{k+1})$.
\item The triple invariant of these new flags
is the reciprocal of the triple invariant of the
original flags.
\end{enumerate}
What this means is that if we have a generic prism, we can
always normalize it so that the corresponding flags are
as in Figure 5.1, with $a_1$ between $a_3$ and $b_1$.
See \S \ref{discuss} for more details.

\subsection{The Big Calculation}
\label{BIG}

   In this section will compute $g^2$,
   the element from Equation \ref{key}, and deduce
   information from the computation.
We first treat the generic case, and then discuss
the non-generic cases at the end of the section.  We normalize as in
Figure 5.1.
\newline
\newline
\noindent
{\bf The Flags:\/}
We represent our points by $3$-vectors in Mathematica:
\begin{equation}
  a_1=(1,0,1) \hskip 30 pt
  a_2=(-1/2,\sqrt 3/2,1), \hskip 30 pt
  a_3=(-1/2,-\sqrt 3/2,1).
\end{equation}
The lines in Figure 5.1 are represented by the cross products
$L_k=a_{k-1} \times a_{k+1}$.
Next, we choose $t>0$ and define
\begin{equation}
  b_k=(1+t) a_k - t a_{k-1}.
\end{equation}
Our flags are $(b_k,L_{k+1})$ for $k=1,2,3$.
The flag fixed by $g^2$ will, as above, be $f_1=(b_1,L_2)$.
The triple invariant of $f_1,f_2,f_3$ is
\begin{equation}
  \label{beep}
  -\frac{t^3}{(t+1)^3} \in (-1,0).
\end{equation}

\noindent
{\bf The Order 3 Element:\/}
The element $\rho(\sigma_3)$ is represented by the matrix
\begin{equation}
M_3=  \left(\matrix{-1/2 & -\sqrt 3/2 & 0 \cr
    \sqrt 3/2 & -1/2 & 0 \cr
    0&0&1}\right)
\end{equation}
This map has order $3$ and has the action $f_k \to f_{k+1}$.
\newline
\newline
{\bf The Order 2 Element:\/}
$\Delta$ be the standard polarity.
Let $S$ be the matrix whose column vectors are
$2rb_1,2sb_2,a_1$:
\begin{equation}
 S= \left(\matrix{
      r(2+3t) & -s (1+3t) &1 \cr
      \sqrt 3 rt & \sqrt 3s(1+t) &0 \cr
      2r & 2s &1}\right)
\end{equation}

Let $L_X$ and $L_Y$ respectively denote
the lines in $\P$ extending the
$X$-axis and the $Y$-axis.  Let
$p_X \in L_X$ and $p_Y \in L_Y$ be the
points at infinity.  The duality
$\Delta$ interchanges the flags
$(p_X,L_X)$ and $(p_Y,L_Y)$ and the
projective transformation represented
by $S$ carries these flags to $f_1$ and $f_2$.
The composition
\begin{equation}
  \label{M2}
  \rho(\sigma_2)=S \circ \Delta \circ S^{-1}= \Delta  \circ M_2 
  \hskip 30 pt 
  M_2=(S^{-1})^t S^{-1}.
  \end{equation}
gives the general form of the elliptic polarity
which interchanges $f_1$ and $f_2$.
Thus, choosing the parameters $(r,s)$ picks
out a generating point in the flat $F$
determined by these two flags.
We compute
\begin{equation}
  \det(S)=6 \sqrt 3 rs t(1+t).
\end{equation}
Since $t>0$, this determinant is nonzero as long as
$rs \not =0$.

Now we observe a symmetry.   If $D$ is any diagonal
matrix whose diagonal entries belong to the
$2$ element set $\{-1,+1\}$ then $D \circ \Delta=\Delta \circ D$.
For this reason, the matrix $S \circ D$ gives the same polarity
as the matrix $S$. This means that all the possibilities
are covered by the cases $r,s>0$.
\newline
\newline
{\bf The Key Element:\/}
Finally, we have
\begin{equation}
  g^2=(M_2^{-1})^t (M_3^{-1})^t M_2 M_3.
  \end{equation}
  To see why this works, we work from right to left.
  We start out with a vector representing a
point.
We apply $M_3$ and we get another vector representing a point.
Now we apply $\Delta M_2$ and we get a vector representing
a line.   The next two matrix operations involve the inverse
transpose because we are acting on lines.  Finally, we apply
$\Delta$ and we get a vector representing a point.

As a sanity check, we compute that
\begin{equation}
  \det(g^2)=1, \hskip 30 pt g^2(b_1)=b_1, \hskip 30 pt
  g^2(L_2)=L_2.
\end{equation}
For the final calculation, of course, we use the inverse
transpose of the matrix representing $g^2$.
Thus $g^2$ fixes the flag $f_1=(b_1,L_2)$, as expected.
\newline
\newline
{\bf Eigenvalues:\/}
Now a miracle occurs.  The matrix $g^2$ is huge, but we compute in
Mathematica that its eigenvalues are:
\begin{equation}
  \label{lambda}
  1, \hskip 10 pt \lambda, \hskip 10 pt \lambda^{-1}, \hskip 30 pt
    \lambda = -\bigg(\frac{r^2}{s^2}\bigg)\bigg( \frac{t}{1+t}\bigg)<0.
\end{equation}
This element is loxodromic unless $\lambda=-1$.
The eigenvector corresponding to $\lambda$ represents $b_1$.
\newline
\newline
{\bf Exploring the Dichotomy:\/}
We have $r,s,t>0$.  From the calculation above, we see that
$g^2$ is loxodromic unless
\begin{equation}
  \label{mu}
  r=\mu s, \hskip 30 pt \mu=\sqrt{\frac{1+t}{t}}.
\end{equation}
We hold $t$ fixed for the rest of this subsection.
With $t$ fixed, the parabolic case is parametrized by the
infinite set $s>0$, which is homeomorphic to a line.
Call this set $\cal J$.

Now let us look at the Pappus representations.
In [{\bf S1\/}] we prove that these representations
correspond to prism pairs $(\Pi,p)$ with $p$ on an
inflection line.  We give a different proof below
in \S \ref{inflectionproof}.
The triple invariant of the representation is an injective
function of our parameter $t$.  With $t$ fixed,
we get an iso-prismatic family parametrized by the
inflection line in $F$: The prisms are all the same.
Call this family $\cal I$.   Each member of $\cal I$
gives us a triple $(r,s,t)$, and this triple must lie
in $\cal J$ because the corresponding $g^2$ is
parabolic. This gives us a map ${\cal I\/} \to {\cal J\/}$.
  No two distinct representations in $\cal I$ are conjugate
  to each other, because a conjugacy would preserve
    the pattern of flats, prisms, triangles, and inflection points.
    Hence, the map
    ${\cal I\/} \to {\cal J\/}$ is injective.
  The family $\cal I$ is parametrized by the point $p$ on the
  inflection line. As $p$ exits every compact subset of $X$,
    the corresponding parabolic element exits every
  compact subset of $SL_3(\R)$.  Hence
  the map ${\cal I\/} \to {\cal J\/}$ is proper.
  Hence  ${\cal J\/}={\cal I\/}$.  This proves
  Theorem \ref{loxo} in the generic case.
    \newline
  \newline
  {\bf The Second Prism Description}
Consider the flag
  $f'_1=(b'_1,L_2')$ corresponding to the eigenvalue
  $1/\lambda$ of $g^2$ and $(g^{-2})^t$.
  The flag $f_1'$ is distinct from $f_1$ because
  the eigenvalues are different.
  The eigenflag $f_1'$ has a monstrous formula,
  but the coordinates are
  rational functions of $r,s,t$.
    We get a new triple of flags by taking the
  orbit of $(p_1',L_2')$ under the action of
  $M_3$.  Thus $p_{k+1}'=M_3(p_k')$ and
  $L_{k+1}=M_3(L_k')$.
  Normallly we would use the inverse-transpose
  to compute the new lines, but in this case
  $(M_3^{-1})^t=M_3$.
  
  The new triple of flags in turn defines
  a new prism $\Pi'$ together with a new
  flat $F'$ of $\Pi'$ corresponding the
  flags $(b_1',L_2')$ and $(b_2',L_3')$.
  We compute (in Mathematica) that $\rho(\sigma_2)$
    swaps $(b_1',L_2')$ with $(b_2',L_3')$.  This
  means that the fixed point of $\rho(\sigma_2)$ in $X$, namely the
  generating
    point $p$ for our representation,
  also lies in $F'$.   So, the pair
  $(\Pi',p)$ is a second description of the same
  prism group.   Exactly one prism pair
  is attracting and one is repelling.
  This proves Theorem \ref{loxo2}
  in the generic case.
    \newline
  \newline
  {\bf The Non-Generic Cases:\/}
    For the totally symmetric case we are back in the
  hyperbolic plane with the Farey triangulation and its
  shears.  In this case,
  Theorem \ref{loxo} follows from the
    hyperbolic geometry picture developed in \S 2.

  The remaining cases correspond to the case
  when the triple invariant is $-1$ but the triple
  is not completely symmetric.  In this case
  we set $b_k=a_k$ for $k=1,2,3$, and
  $L_1,L_2,L_3$ are the line through the origin
  with $b_k \in L_{k+1}$.
  The matrix $S$ above is now the one
  whose column vectors are $2rb_1,2rb_2,(0,0,1)$.
  With these changes, the calculation above, and
  all the results, go through just as in the generic
  case.
  
  \subsection{Comparing the Prisms}
  \label{comparexx}

    We consider the loxodromic case in more detail.
        We call the two prism pairs $(\Pi,p)$ and
    $(\Pi',p)$ {\it partners\/}.  
  Equation \ref{beep} gives the triple
  invariant for the flags defining $\Pi$.
  The invariant for $\Pi$
  is $3\log((t+1)/t)$.
    Let $\tau'$ be the triple invariant for the flags
    defining $\Pi'$.  See \S \ref{monster} below for the
        monstrous formula.   Inspecting this formula,
    we see that
    $\tau'<0$ no matter which $r,s,t>0$ we choose.
    See the discussion at the end of \S \ref{monster}.
    As in Equation \ref{beep},
    the expression for $\tau'$ is a perfect cube.

Here is a sample calculation.
  The prism invariants for $\Pi$ and  $\Pi'$ when $(r,s,t)=(1,1,1)$ are
  respectively
  $$
  3 \log(2), \hskip 30 pt
  3\log \bigg(\frac{9825}{5602}\bigg)
  $$

  In the non-generic case, the formula for $\tau'$ is shorter:
  {\tiny
$$
-\frac{\left(36 r^{12} s^2+r^{12}+3 r^{10} s^2+3 r^8 s^4+1296 r^6
    s^{10}+144 r^6 s^8+2 r^6 s^6+108 r^4 s^{10}+3 r^4 s^8+3 r^2
    s^{10}+s^{12}\right)^3}{\left(r^{12}+1296 r^{10} s^6+108 r^{10}
    s^4+3 r^{10} s^2+144 r^8 s^6+3 r^8 s^4+2 r^6 s^6+3 r^4 s^8+36 r^2
    s^{12}+3 r^2 s^{10}+s^{12}\right)^3}$$
\/}
One can see, again, that this expression is always negative.
When $r=s$ the expression equals $-1$ as it must in the
parabolic case.  In general, the expression can take on all
negative values.  Reversing the roles played by
$\Pi$ and $\Pi'$ we see that a non-generic pair
$(\Pi',p)$ can arise from a generic pair $(\Pi,p)$ no
matter what the prism invariant of $\Pi$.  It all
depends on the choice of $p$.

\subsection{Pappus Representations and the Inflection Line}
\label{inflectionproof}

In this section we prove, as we did in [{\bf S1\/}],
that the Pappus representations correspond
to prism pairs $(\Pi,p)$ where $p$ lies on the inflection line of the
flat $F$ which contains $p$.  Our proof is a bit different here.
In \S \ref{AXIAL} we prove:

\begin{lemma}
  \label{INFL}
  For each prism $\Pi$ there exists a Pappus modular group
  whose order  $2$ fixed points are the inflection points of $\Pi$.
\end{lemma}

Recall that the flat $F$ is foliated
by medial geodesics defined by the flags
$(b_1,L_2)$ and $(b_2,L_3)$.   The orthogonal
foliation is by singular geodesics which we call
{\it orthgonal singular geodesics\/} .

\begin{lemma}
  \label{ortho}
  \label{ortho2}
  The following is true.
  \begin{enumerate}
    \item
  The points $(r,s,t)$ and $(rd,sd,t)$
  parametrize representations which
  are $\sqrt{2/3} \log(d)$ units apart along
  an orthogonal singular geodesic.
\item
  The points $(r,s,t)$ and $(rd,s/d,t)$
  parametrize representations which
  are $\sqrt 2\log(d)$ units apart along
  a medial geodesic.
  \end{enumerate}
\end{lemma}

\startproof
Referring to Equation \ref{M2}, let
$I_0=M_2(r,s,t)$ and
$I_1=M_2(rd,sd,t)$.  Let
$$J=(I_1^{-1})^tI_0.$$
Then $J$ is a translation along
$F$ along the vector which is twice
the difference between the fixed
points of $I_0$ and $I_1$.
We compute that the eigenvalues
of $J$ are
$(1,d^2,d^2)$, and
the vectors representing $b_1$ and $b_2$
are both eigenvectors corresponding to
$d^2$.   This shows that $J$ is translation
along the singular geodesics by
$\sqrt{8/3} \log d$ units.
(We get the funny factor because we compute
the distance by first rescaling $J$ so that it
belongs to $SL_3(\R)$.)
Hence the fixed points of $I_0$ and $I_1$ are
$\sqrt{2/3}\log(d)$ units apart and on the same orthogonal
singular geodesic.  This proves Statement 1.

For Statement 2, we set
$I_1=M_2(rd,s/d,t)$.
We compute that the eigenvalues
of $J \in SL_3(\R)$ are
$(d^2,d^{-2},1)$, and
the vectors representing $b_1$ and $b_2$
are eigenvectors respectively corresponding
to $d^2$ and $1/d^2$.  Finally, the
point $a_1=L_2 \cap L_3$ is the eigenvector
corresponding to $1$.   Our result follows
from this structure.
\endproof

Our claim about the inflection line follows from
 Equation \ref{mu}, Lemma
 \ref{INFL}, and Statement 1 of Lemma \ref{ortho}.
 Statement 2 of Lemma \ref{ortho} is not needed until later.

\subsection{Symmetries of the Parametrization}
\label{discuss}

Now we discuss some of the symmetries of our $(r,s,t)$
parametrization.
\newline
\newline
  {\bf First Symmetry:\/}
Define
  \begin{equation}
    \iota_1(r,s,t)=\bigg(\frac{r_0^2}{r},\frac{s_0^2}{s},t\bigg).
  \end{equation}
  Here $r_0$ and $s_0$ depend on $t$ but we
  have suppressed the dependence.
  Note that the map $\iota_1$ is an involution
   which fixes $(r_0,s_0,t)$.  If we
  use $\log$-coordinates and identify
  the $(r,s)$ positive quadrant with
  $\R^2$, then this map is the order $2$
  isometry  of $\R^2$ which fixes the point
  corresponding to the axial representation.
  This map corresponds to the action of the
  order $2$ symmetry of the prism $\Pi$ which
  preserves the flat $F \cong \R^2$ containing
  the generating point. The representations
  corresponding to $(r,s,t)$ and $\iota_1(r,s,t)$ are
  not conjugate, but they have the same
  image in ${\rm Isom\/}(X)$.   The representations
    become conjugate if we switch the
  order $3$ generator of one of the
  groups to its inverse. \newline
  \newline
{\bf Second Symmetry:\/}
In \S \ref{normal} we said that
we would restrict our attention to
  triples of flags with triple invariant in $(-1,0)$.  This
  corresponds to $t \in (0,\infty)$.  The case of
  flags with triple invariant in $(-\infty,-1)$ corresponds
  to $t \in (-\infty,-1)$.   We mentioned that we can
  always arrange the former case by applying a duality
  to our initial triple of flags. 
  We define
  \begin{equation}
    \label{DI}
    \iota_2(r,s,t)=\bigg(\frac{s_0^2}{s},\frac{r_0^2}{r},-1-t\bigg).
  \end{equation}
  The two sets of parameters $(r,s,t)$ and
  $\iota_2(r,s,t)$ describe groups that are conjugate
  {\it via\/} a duality.
  We call $\iota_2$ 
  the {\it duality involution\/}.

\subsubsection{The Monster Expression}
\label{monster}

The expression for the second triple invariant $\tau'$ is
$-(1+t)^3 A^3/(t^3 B^3)$, where $A$ and $B$ respectively are:

{\footnotesize
$$
 \matrix{
576 r^{12} s^2 t^8 + 1152 r^{12} s^2 t^7 + 960 r^{12} s^2 t^6 + 384 r^{12} s^2 t^5 + 64 r^{12} s^2 t^4 + 16 r^{12} t^8 + \cr
48 r^{10} s^2 t^8 + 192 r^{10} s^2 t^7 + 288 r^{10} s^2 t^6 + 192 r^{10} s^2 t^5 + 48 r^{10} s^2 t^4 + 48 r^8 s^4 t^8 + \cr
96 r^8 s^4 t^7 + 48 r^8 s^4 t^6 + 12 r^8 s^2 t^6 + 24 r^8 s^2 t^5 + 12 r^8 s^2 t^4 + 20736 r^6 s^{10} t^8 + \cr
82944 r^6 s^{10} t^7 + 152064 r^6 s^{10} t^6 + 165888 r^6 s^{10} t^5 + 117504 r^6 s^{10} t^4 + 55296 r^6 s^{10} t^3 + \cr
16896 r^6 s^{10} t^2 + 3072 r^6 s^{10} t + 256 r^6 s^{10} + 2304 r^6 s^8 t^8 + 4608 r^6 s^8 t^7 + 3840 r^6 s^8 t^6 + \cr
1536 r^6 s^8 t^5 + 256 r^6 s^8 t^4 + 32 r^6 s^6 t^8 {\bf - 192 r^6 s^6 t^7 - 192 r^6 s^6 t^6 - 64 r^6 s^6 t^5\/} + 16 r^6 s^4 t^6 + \cr
r^6 s^2 t^4 + 1728 r^4 s^{10} t^8 + 10368 r^4 s^{10} t^7 + 27072 r^4 s^{10} t^6 + 40320 r^4 s^{10} t^5 + 37632 r^4 s^{10} t^4 + \cr
22656 r^4 s^{10} t^3 + 8640 r^4 s^{10} t^2 + 1920 r^4 s^{10} t + 192 r^4 s^{10} + 48 r^4 s^8 t^8 + 192 r^4 s^8 t^7 + 288 r^4 s^8 t^6 + \cr
192 r^4 s^8 t^5 + 48 r^4 s^8 t^4 + 48 r^2 s^{10} t^8 + 384 r^2 s^{10} t^7 + 1344 r^2 s^{10} t^6 + 2688 r^2 s^{10} t^5 + \cr
3360 r^2 s^{10} t^4 + 2688 r^2 s^{10} t^3 + 1344 r^2 s^{10} t^2 + 384 r^2 s^{10} t + 48 r^2 s^{10} + 16 s^{12} t^8 + \cr
96 s^{12} t^7 + 240 s^{12} t^6 + 320 s^{12} t^5 + 240 s^{12} t^4 + 96 s^{12} t^3 + 16 s^{12} t^2 + 4 s^{10} t^6 + \cr
24 s^{10} t^5 + 60 s^{10} t^4 + 80 s^{10} t^3 + 60 s^{10} t^2 + 24 s^{10} t + 4 s^{10}
}
$$

$$
\matrix{
r^2 s^6 + 16 r^4 s^6 + 96 r^6 s^6 + 256 r^8 s^6 + 256 r^{10} s^6 + 16 s^{12} + 64 r^2 s^{12} + \cr
4 r^2 s^6 t + 96 r^4 s^6 t + 768 r^6 s^6 t + 2560 r^8 s^6 t + 3072 r^{10} s^6 t + 128 s^{12} t + \cr
640 r^2 s^{12} t + 6 r^2 s^6 t^2 + 240 r^4 s^6 t^2 + 2688 r^6 s^6 t^2 + 11520 r^8 s^6 t^2 + 16896 r^{10} s^6 t^2 + \cr
12 r^2 s^8 t^2 + 48 r^4 s^8 t^2 + 448 s^{12} t^2 + 2880 r^2 s^{12} t^2 + 4 r^2 s^6 t^3 + 320 r^4 s^6 t^3 + 5312 r^6 s^6 t^3 + \cr
30208 r^8 s^6 t^3 + 55296 r^{10} s^6 t^3 + 48 r^2 s^8 t^3 + 288 r^4 s^8 t^3 + 896 s^{12} t^3 + 7552 r^2 s^{12} t^3 + \cr
48 r^8 s^4 t^4 + 192 r^{10} s^4 t^4 + r^2 s^6 t^4 + 240 r^4 s^6 t^4 + 6400 r^6 s^6 t^4 + 50176 r^8 s^6 t^4 + \cr
117504 r^{10} s^6 t^4 + 72 r^2 s^8 t^4 + 720 r^4 s^8 t^4 + 48 r^2 s^{10} t^4 + 1120 s^{12} t^4 + 12544 r^2 s^{12} t^4 + \cr
192 r^8 s^4 t^5 + 1152 r^{10} s^4 t^5 + 96 r^4 s^6 t^5 + 4736 r^6 s^6 t^5 + 53760 r^8 s^6 t^5 + 165888 r^{10} s^6 t^5 + \cr
48 r^2 s^8 t^5 + 960 r^4 s^8 t^5 + 192 r^2 s^{10} t^5 + 896 s^{12} t^5 + 13440 r^2 s^{12} t^5 + 4 r^{10} t^6 + \cr
16 r^{12} t^6 + 288 r^8 s^4 t^6 + 2880 r^{10} s^4 t^6 + 16 r^4 s^6 t^6 + 2048 r^6 s^6 t^6 + 36096 r^8 s^6 t^6 + \cr
152064 r^{10} s^6 t^6 + 12 r^2 s^8 t^6 + 720 r^4 s^8 t^6 + 288 r^2 s^{10} t^6 + 448 s^{12} t^6 + 9024 r^2 s^{12} t^6 + \cr
32 r^{12} t^7 + 192 r^8 s^4 t^7 + 3456 r^{10} s^4 t^7 + 448 r^6 s^6 t^7 + 13824 r^8 s^6 t^7 + 82944 r^{10} s^6 t^7 + \cr
288 r^4 s^8 t^7 + 192 r^2 s^{10} t^7 + 128 s^{12} t^7 + 3456 r^2 s^{12} t^7 + 16 r^{12} t^8 + 48 r^{10} s^2 t^8 + \cr
48 r^8 s^4 t^8 + 1728 r^{10} s^4 t^8 + 32 r^6 s^6 t^8 + 2304 r^8 s^6 t^8 + 20736 r^{10} s^6 t^8 + 48 r^4 s^8 t^8 + \cr
48 r^2 s^{10} t^8 + 16 s^{12} t^8 + 576 r^2 s^{12} t^8
}
$$
}

\noindent
{\bf Remarks:\/}
\newline
(1) I used chatGPT 4o mini to help me format these
monsters. \newline
(2)
We have $B>0$ for all $r,s,t$ because $B$ is the sum of positive monomials.
When $r \geq 1$ and $t \geq 1$ and $s \leq 1$, the first term of $A$
dominates the sum of the bolded negative terms, and so $A>0$ in this
case.   There are $7$ other cases, depending on the signs of
$r-1, s-1, t-1$, and in all cases one can find short sums of positive
terms which dominate
the $3$ bolded terms.  Hence $A>0$ for all $r,s,t>0$.
Hence $\tau'<0$ for all $r,s,t>0$.

      \section{Recognizing the Representations}

\subsection{The Pappus Modular Groups}

In this chapter we relate the prism representations
to the Pappus modular groups and to the Anosov
representations considered in
[{\bf BLV\/}] and [{\bf S2\/}].   This section is a
near-repeat of the corresponding section in
[{\bf S2\/}].   One subtle difference is that
we choose the parameters $c,d$ in $(-1,1)$
as in [{\bf S2\/}] rather than in $(0,1)$ as in
[{\bf S1\/}].   
\newline
\newline
{\bf Convex Marked Boxes:\/}
A {\it convex marked box\/} is a convex quadrilateral
in $\P$ together with a distinguished point in the
interior of one side and a distinguished point in the
interior of an opposite side.  We call one of the
points the {\it top\/} point and the other one the
{\it bottom\/} point.  Correspondingly we call the
edges containing these points the {\it top edge\/}
and the {\it bottom edge\/}.  Finally, we say that
the {\it top flag\/} is the flag $(p,\ell)$ where $p$
is the top point and $\ell$ is the line extending the
top edge.  We define the {\it bottom flag\/} similarity.
\newline
\newline
{\bf Operations on Marked Boxes:\/}
There are $3$ operations we can perform on marked
boxes, and we call them $t,b,i$.  Figure 6.1
shows how they act.

 \begin{center}
\resizebox{!}{1.7in}{\includegraphics{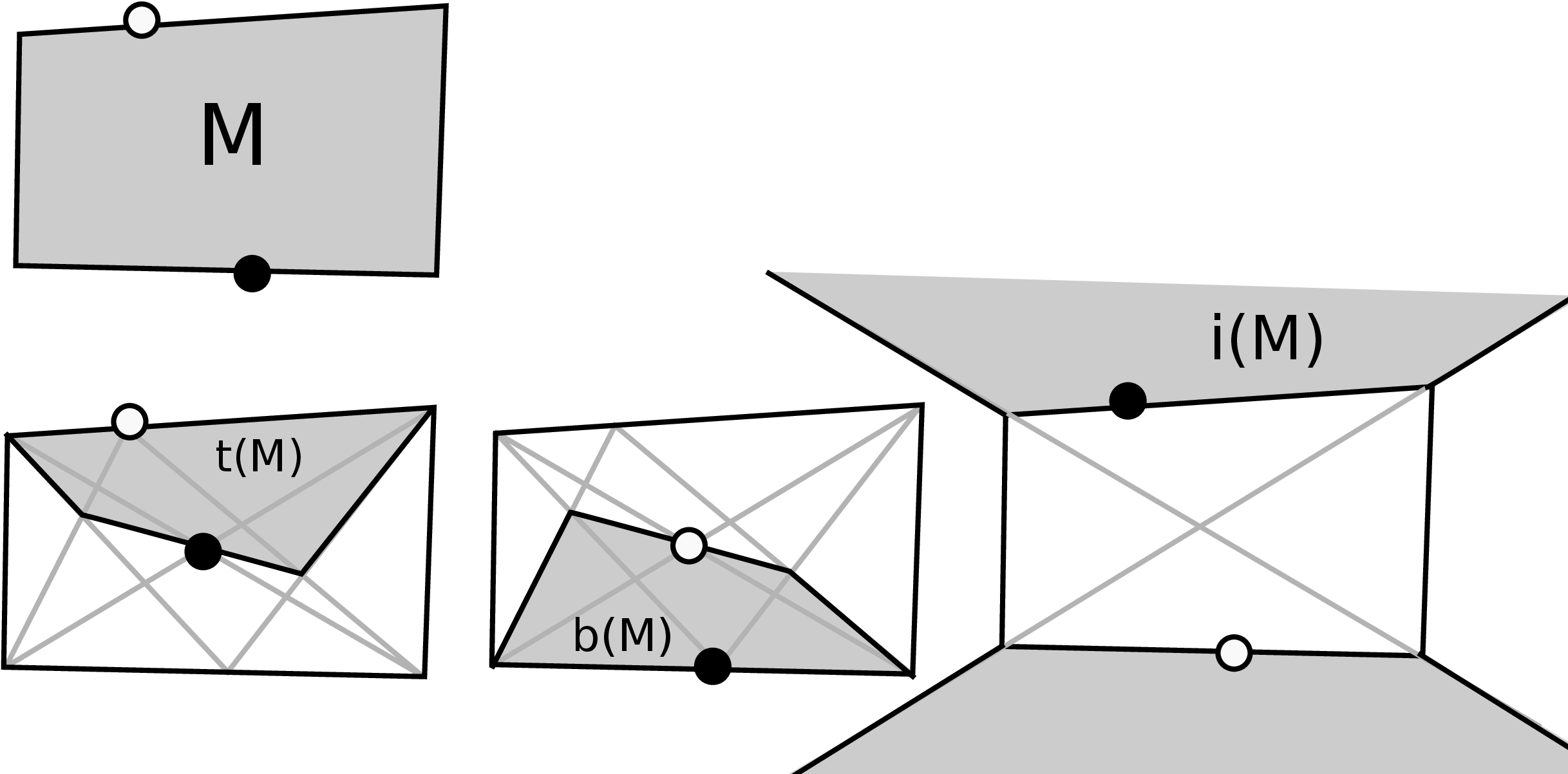}}
\newline
{\bf Figure 6.1:\/} The three operations on marked boxes
packing
\end{center}

These operations satisfy the relations
\begin{equation}
  i^2=I. \hskip 20 pt
  tit=b, \hskip 20 pt bib=t, \hskip 20 pt
  tibi=I, \hskip 20 pt biti=I.
\end{equation}
here $I$ is the identity.  As a consequence of these relations,
and the nesting of the marked boxes.  The group of operations
isomorphic to the modular group.    The explicit generators are
(say) $i$ and $ti$.  
We let $\cal M$ be the orbit of a marked box
under the action of this group.
\newline
\newline
{\bf Order Three Symmetries of the Orbit:\/}
Given a marked box $M \in \cal M$ there is an order $3$ projective
transformation $T_M$ which has the orbit
$$i(M) \to t(M) \to b(M).$$
This accounts for the order
$3$ elements of the Pappus modular groups.  If
we list out the top and bottom flats of these three marked
boxes, they coincide in pairs and we end up with a
triple of flags.  The triple always turns out to be negative.
Thus each marked box $M$ in $\cal M$ gives us a prism
in $X$.
\newline
\newline
{\bf Order Two  Symmetries of the Orbit:\/}
There is also an elliptic polarity which, in a certain sense,
swaps $M$ and $i(M)$.  To make sense of this, we have
to recall the notion of a {\it doppelganger\/} defined in
[{\bf S1\/}].

 \begin{center}
\resizebox{!}{1.2in}{\includegraphics{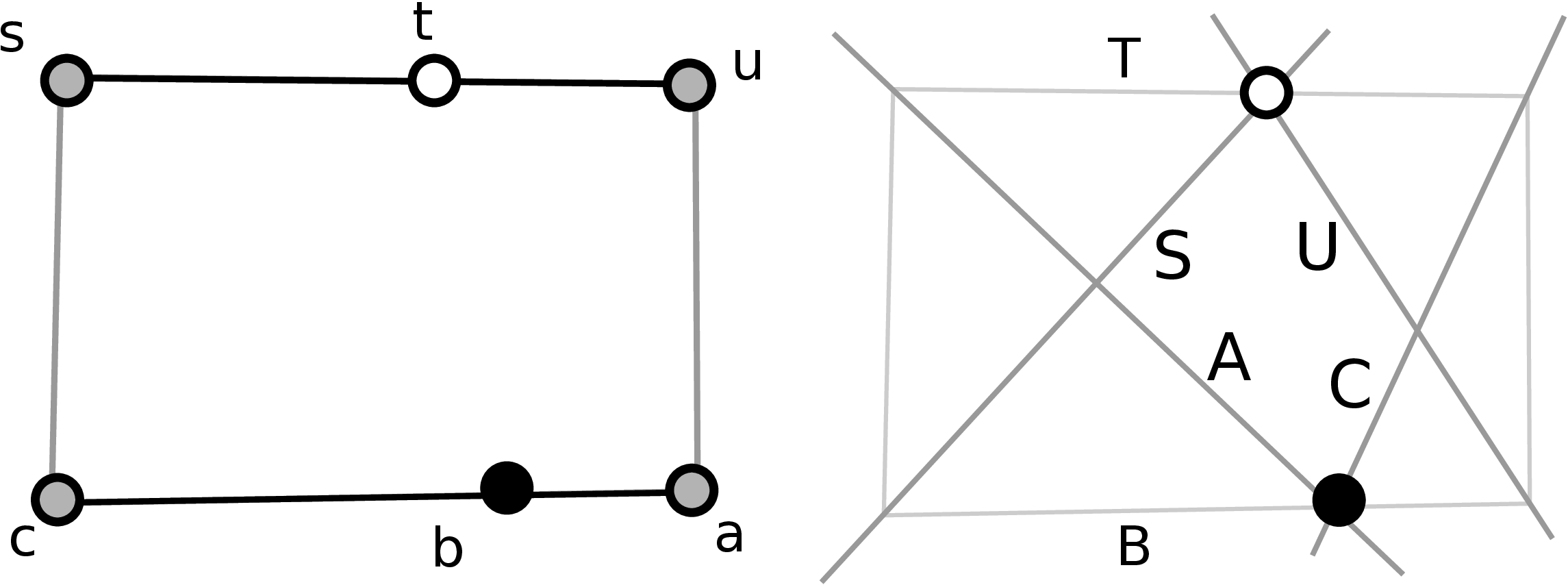}}
\newline
{\bf Figure 6.2\/} A convex marked box and its doppelganger
\end{center}

The $6$-tuple $(s,t,u,a,b,c)$ shown on the left side of Figure 6.2
encodes the marked box $M$.  Here $t$ and $b$ are respectively
the top and bottom points of $M$.  The corresponding $6$-tuple of lines
$(S,T,U,A,B,C)$, which is defined entirely in terms of $M$,
encodes a convex marked box $M^*$ in $\P^*$.    We can repeat
the operation and we get $M^{**}=M$.    It turns out that the
$i,b,t$ operations commute with the doppelganger operation and
we can think of our orbit $\cal M$ as an orbit of pairs of the
form $(M,M^*)$.  We call such a pair an {\it enhanced convex marked
  box\/}.

We showed in [{\bf S1\/}] that there is an elliptic polarity
$\delta_M$ 
that swaps $M$ and $(i(M))^*$, and simultaneously
swaps $M^*$ and $i(M)$.
\newline
\newline
{\bf Putting it Together:\/}
Given a marked box $m$ we have
$3$ flags coming from the tops and bottoms of
$i(m)$, $t(m)$ and $b(m)$.
We let $\Pi$ be the corresponding prism.
To get a prism representation $\rho$, 
we let $\rho(\sigma_2)=\delta_M$ and we
let $\rho(\sigma_3)$ be the order $3$ symmetry
discussed above.  By construction
$\rho(\sigma_3)$ permutes the flats of $\Pi$, and
the fixed point $p$ of $\rho(\sigma_2)$ lies in
flat $F$ of $\Pi$ determined by the tops and
bottoms of $i(m)$.

\subsection{The Axial Case}
\label{AXIAL}

We define an equivalence relation on real
numbers, $(c,d) \sim (-c,-d)$.
We denote the equivalence class of
$(c,d)$ by $[(c,d)]$.  A marked box
$M$ has an invariant
$[M]=[(c,d)]$ where $c$ and $d$ are
defined as follows.   We normalize by
a projective transformation so that
the quadrilateral underlying $M$ has
vertices $(\pm 1,\pm 1)$.
The top point then has coordinates $(c,1)$
and the bottom point has coordinates $(d,-1)$.
The reason for the equivalence relation is that
we might have normalized by post-composing
our projective transformation by the map
which is reflection in the vertical line $x=0$.

We say that $M$ is {\it axial\/} if
$[M]=[(0,d)]$ or $[M]=[(c,0)]$.
Here is a sharper
version of Lemma \ref{INFL}.  We also
prove this result in [{\bf S1\/}].

\begin{lemma}
  Let $M$ be an axial marked box and let
  $(\Pi,p)$ be the prism pair corresponding to the
  associated Pappus modular group representation.
  Then $p$ is an inflection point of $\Pi$.
\end{lemma}

\startproof
We treat the case when $[M]=[(0,d)]$. The
other case follows from duality: $[i(M)]=[(d,0)]$.

Two marked boxes $M_1$ and $M_2$
are projectively equivalent if and only if
$[M_1]=[M_2]$.   Also, if $[M]=[(c,d)]$ then
$i(M)=[(-d,c)]$.   This is a direct
calculation \footnote{Again we mention that we
  normalized so that $(c,d) \in (0,1)$ in [{\bf S1\/}].
  Here we are taking $(c,d) \in (-1,1)$ and so the
  formulas are a bit different.}
which we do in [{\bf S1\/}]. So, if
$[M]=[(0,d)]$ then, as we already mentioned,
$i(M)=[(-d,0)]=[(d,0)]$.
Hence there is a projective transformation $\phi$
which maps $M$ to $i(M)$, switching the tops and
bottoms.   Note that $\phi^2$ fixes all the vertices
of $M$, and hence $\phi^2$ is the identity.

For notational convenience set $h=\rho(\sigma_2)$.
Define $M'$ to be the marked box we get by
simply exchanging what we mean by the top and
the bottom.  We let
$\psi=h\circ \phi$.  Then
$\psi$ is a duality which maps $M$ to
the doppelganger of $M'$.
So, $\psi$ maps $t(M)$ and $b(M)$ respectively
to the doppelgangers of $b(M)$ and $t(M)$.
In particular, $\psi$ permutes the $3$ flags
associated to $i(M), t(M), b(M)$. Hence
$\psi(\Pi)=\Pi$.   We have identified $\psi$
as symmetry of $\Pi$.  Since $\psi$ is a duality,
$\psi$ has order $2$ and fixes one of the
inflecton points of $\Pi$.

Since $h=h^{-1}$ we have $\phi= h\psi$.
But then $h\psi h\psi$ is the identity.
Hence $h\psi=\psi h$.
Hence $\psi$ and $h$ have the same fixed point.
Since the fixed point of $\psi$ is an inflection
point of $\Pi$, the fixed point of $h$ is the
same infection point.
\endproof

\subsection{The Space of Prism Representations}

Two Pappus modular group representations are conjugate if and
only if the enhanced marked boxes in their orbits are
projectively equivalent, either by dualities or projective
transformations.   We can get a section of the space of
Pappus representations by normalizing so that
our initial marked box is the square $Q$
with vertices $\pm 1$ and the top
point $t$ lies in the interior of the top edge of $Q$ and
the bottom point $b$ lies in the interior of the bottom
edge of $Q$.

We let $M(c,d)$ be the marked box which is normalized
this way and has
$[M]=[(c,d)]$.
The boxes $M(c,d)$ and $M(c',d')$ define the same
representation in $\cal R$ if and only if
$(c,d)$ and $(c',d')$ are in the same $\rho$-orbit,
where $\rho$ is the order $4$ rotation about the origin.
This, as we saw in
[{\bf S1\/}], the space of Pappus modular group
representations is homeomorphic to the cone
\begin{equation}
{\cal C\/}=(-1,1)^2/\langle \rho \rangle.
\end{equation}
This space is in turn homeomorphic to $\R^2$.

Let $h$ be the map which assigns to each point
in $\cal C$ the isometry class of prism pairs
for the associated Pappus modular group representation.
The map $h$ is a homeomorphism between $\cal C$ and $\cal P$.

\begin{lemma}
The map $h$ extends to be a homeomorphism
from ${\cal C\/} \times [0,\infty)$ to $ {\cal P\/} \cup {\cal A\/}_+$.
\end{lemma}

\startproof
Given $c \in \cal C$ and $r \in [0,\infty)$
let $(\Pi,p_0)$ be the prism pair given by $h(c)$. 
Let $F$ be the flat of $\Pi$ containing $p_0$.
Recall that $p_0$ lies on the inflection line in $F$.
Let $\gamma$ be the geodesic in the medial geodesic foliation
of $F$ that contains $p_0$.  One component of
$\gamma-p_0$ consists of points $q$ such that
$(\Pi,q)$ is attracting and the other component
is the repelling case.  These components cannot mix
because the only neutral pairs lie on the inflection line.
So, we let $q_r$ be the unique point of $\gamma$ such that
$d_X(p_0,p_r)=r$ and $(\Pi,p_r)$ is an attracting pair.

By construction, $h$ gives a continuous proper bijection from
${\cal C\/} \times [0,\infty)$ to ${\cal P\/} \cup {\cal A\/}_+$.
A map with all these properties is a homeomorphism.
The continuity follows from the fact that if
$(\Pi,p_0)$ and $(\Pi',p_0')$ are two nearby prism
pairs, then the points $p_r$ and $p_r'$ are also
close in $X$.   The key observation is that
the attracting rays of $\gamma$ and $\gamma'$
point in the about the same rather than about
opposite directions.   The properness follows from
the fact that, as $r \to \infty$, the distance from
$p_r$ to the inflection line of $\Pi$ tends to $\infty$.
\endproof

\subsection{The Image in the Big Representation Space}
\label{IMAGE}

Let $O=(0,0,0)$, the origin in $\R^3$.
We have a composition of maps
\begin{equation}
  \R^2 \times [0,\infty) \to  {\cal P\/} \cup {\cal A\/}_+ \to {\cal
    R\/} \cong \R^3-O.
\end{equation}
The first of these maps is a homeomorphism.
The second of these maps is both continuous and
injective.  Therefore, the composition
\begin{equation}
  f: \R^2 \times [0,\infty) \to \R^3-O
\end{equation}
is continuous and injective.
 
Since $f$ is a map between $3$-manifolds, it follows
from Invariance of Domain that
$f(\R^2 \times (0,\infty))$ is an open subset of
$\R^3-O$ and $f$ is a homeomorphism from this
set onto its image.  We also remark that the image of $f$
stays outside a neighborhood of $O$
because representations indexed by points close to $O$ satisfy
${\rm trace\/}(g^2) \sim 0$, and all prism representations
have ${\rm trace\/}(g^2) \leq -1$.  Here
$g^2$ is the element we have studied extensively
in the previous chapter.

\begin{lemma}
  $f$ is a proper map.
\end{lemma}

\startproof
What we mean is that if $\{q_n\}$ is a sequence of
points in $\R^2 \times [0,\infty)$ that exits every compact
subset, then $\{f(q_n)\}$ also exits every compact subset
of $\R^3-O$.  Since our image avoids a neighborhood of $O$
we are really saying the the image sequence exits every
compact subset of $\R^3$.  We suppose not and
derive a contradiction.

Let $(\Pi_n,p_n)$ be the prism pair associated to $q_n$.
Let $\eta(\Pi_n,p_n)$ be the invariant computed in
\S \ref{axis}.  The Properness Theorem tells us that
if we have
$\eta(\Pi_n,p_n) \to \infty$ then we also have
$\nu(\Pi_n,p_n) \to \infty$. But this
latter quantity is the distance in $X$ from $p_n$,
the fixed point of the element $\rho_n(\sigma_2)$, to
the geodesic fixed by the element $\rho_n(\sigma_3)$.
Here we are setting $\rho_n=\rho(\Pi_n,p_n)$.
If this distance tends to $\infty$ then our representations
exit every compact subset of $\cal R$.
We conclude that
$\{\eta(\Pi_n,p_n)\}$ remains uniformly bounded.

We want to see that in this case we also have
$\nu(\Pi_n,p_n) \to \infty$.  Since
$\{q_n\}$ is exiting every compact subset of
$\R^2 \times [0,\infty)$ it means that the first
two coordinates of $q_n$ are exiting every
compact subset of $\R^2$.   The corresponding
Pappus modular groups are exiting every compact
subset of $\cal P$.  Since there is a uniform bound
between $p_n$ and the point on the inflection line
contained in the same medial geodesic, it suffices
to prove our result when $\rho(\Pi_n,p_n)$ is a
Pappus representation.   We proved this result
already in [{\bf S2\/}, \S 4].
\endproof

The image of
${\cal A\/}$ in $\cal B$ lies on the side
of $\cal P$ which has no elliptic elements.
Since our map $f$ is proper and injective,
the image of $\cal A$ in $\cal R$ is
precisely $\cal B$.

\section{Patterns of Geodesics and Shearing}

In this chapter we prove
Theorems \ref{extra} and \ref{three}.

\subsection{The Main Argument}

Given the result at the end of the last chapter,
we have $3$ different descriptions of
the Barbot component $\cal B$.
\begin{enumerate}
\setcounter{enumi}{-1}
\item $\cal B$ is the component of representations
  produced by the box morphing construction
  in [{\bf BLV\/}] when it is extended as in [{\bf S2\/}].
From this description, every element of $\cal B$ is
either Pappus or Anosov.

\item $\cal B={\cal A\/}_+ \cup {\cal P\/}$.
  We can make this identification because the map
  from ${\cal A\/} \cup {\cal P\/}$ is a homeomorphism.
  In particular, every element of $\cal B$ is a
  prism group.

\item $\cal B={\cal A\/}_- \cup {\cal P\/}$.
  We can make this identification because the map
  from ${\cal A\/} \cup {\cal P\/}$ is a homeomorphism.
  In particular, every element of $\cal B$ is a
  prism group in a different way.
\end{enumerate}

The space ${\cal P\/} \cup {\cal A\/}_+$ has a
foliation by rays, coming from the medial geodesics which
naturally foliate the prisms.
Indentification 1 transfers this foliation
to a foliation of $\cal B$ by rays.
The space ${\cal P\/} \cup {\cal A\/}_-$ also has a
foliation by rays.  Identificaton 2
transfers this foliation
to a second foliation of
$\cal B$ by rays. 
Our analysis of the triple invariants shows that
generically these two foliations of $\cal B$ do not
coincide.

Let $\rho \in \cal B$.
When we use Identification 1, we produce
prism description $(\Pi_1,p_1)$ of $\rho$.
When we use Identification 2, we produce a second
prism description $(\Pi_2,p_2)$ of $\rho$.
Our analysis of the triple invariants shows that
generically $\Pi_1$ and $\Pi_2$ are not isometric
to each other.  When $\rho \in \cal P$, the
two descriptions coincide.

The prism $\Pi_k$ contains a
distinguished triangle $\gamma_k$,
namely the one that contains the
generating point $p_k$.  The orbit of
$\gamma_k$ under $\Gamma$ is
the pattern $Y_{\rho,k}$.
Here we have set $\Gamma=\rho(\Z/2*\Z/3)$,
as in Theorems \ref{extra} and \ref{three}.
By construction,
$$\Gamma \subset {\rm Isom\/}(Y_{\rho,1}), \hskip 30 pt
  \Gamma \subset {\rm Isom\/}(Y_{\rho,1}).$$
The group $\Gamma$ acts transitively on
the prisms associated to $Y_{\rho,k}$.
Hence each coset of
$\Gamma$ in ${\rm Isom\/}(Y_{\rho,k)}$
has a representative which is a symmetry
of $\Pi_k$.  This shows that our containment
of groups has index at most $6$.  Generically,
the symmetries of $\Pi_k$ do not preserve
the set of inflection points of $\Pi_k$.
For this reason,
$\Gamma={\rm Isom\/}(Y_{\rho,k})$ in the
generic case.

Since $\Pi_1$ and $\Pi_2$ are not generically
isometric, neither are the patterns
$Y_{\rho,1}$ and $Y_{\rho,2}$.
We also have the patterns
$F_{\rho,1}$ and $F_{\rho,2}$ which
are the Farey patterns corresponding
to the endpoints of the two rays
containing our representation.
The same triple invariant analysis shows that
$F_{\rho,1}$ and $F_{\rho,2}$ are not isometric to each other.
We have not yet pinned down the relationship
between $F_{\rho,k}$ and $Y_{\rho,k}$.  We will
show below that the relationship is
that of shearing.

To finish the proof we need to show that the
patterns $Y_{\rho,k}$ are embedded for $k=1,2$,
and we need to explain the shearing relationship
between $F_{\rho,k}$ and $Y_{\rho_k}$. Finally, we need to
establish the equality of the shearing strengths.
The rest of the chapter is devoted to these things.

\subsection{Pairs of Flags}

In this section we give some
preliminary information about certain
pairs of flags.

We say that a pair of flags is {\it orthogonal\/}
if the flat it determines contains the origin.
In this case, the standard polarity switches
the two flags.  Figure 7.1 shows a typical
orthogonal pair, $(p_1,L_1)$ and $(p_2,L_2)$.
The lines $L_1$ and $L_2$ are parallel, the
the line $\overline{p_1p_2}$ contains the
origin and is perpendicular to $L_1$ and $L_2$.
Finally, we have $d_1d_2=1$ where $d_k$ is the
distance from $L_k$ to the origin.   

   \begin{center}
     \resizebox{!}{2.2in}{\includegraphics{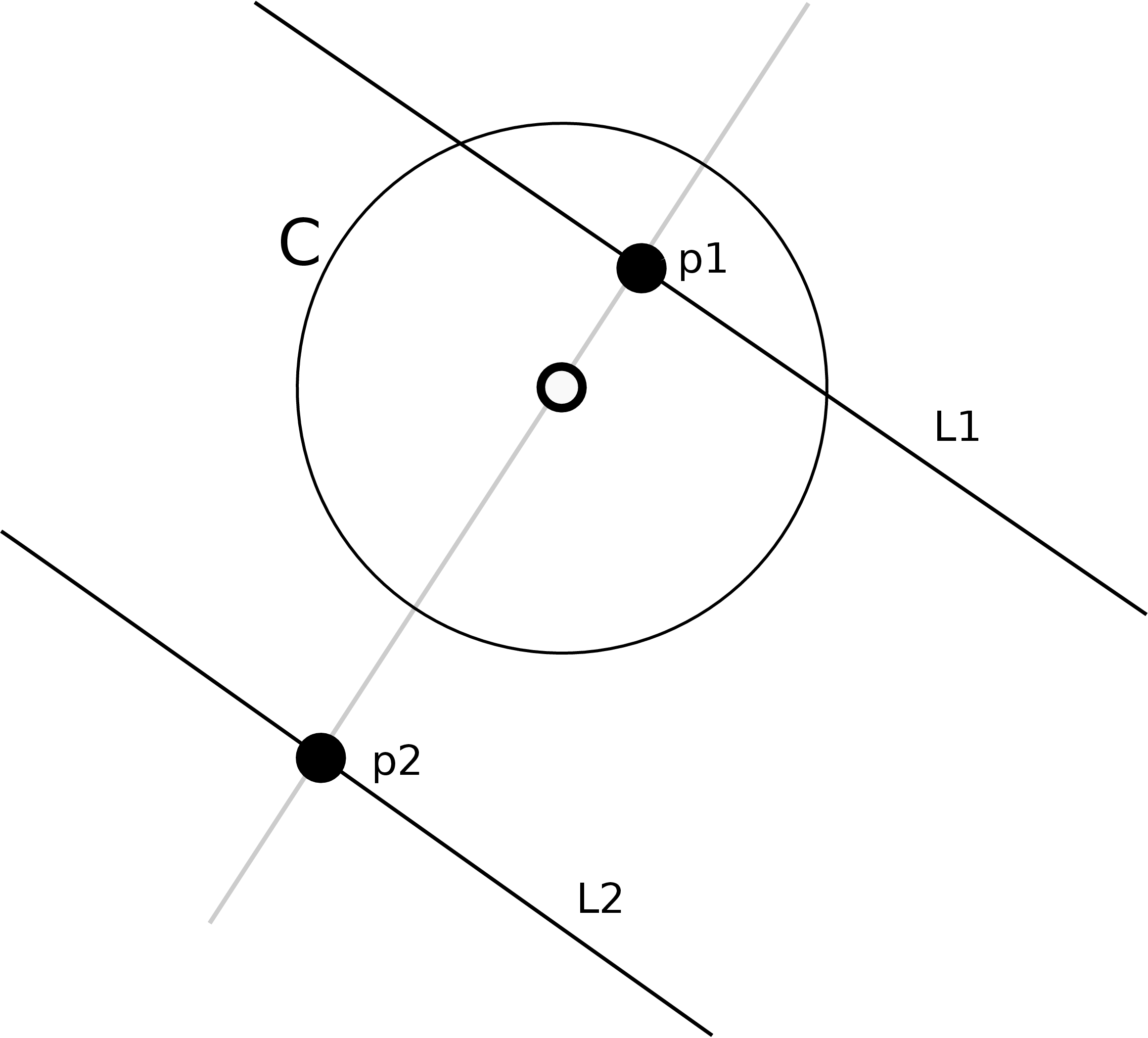}}
    \newline
    {\bf Figure 7.1\/}  An orthogonal pair of flags and the unit
    circle $C$.
      \end{center}

      If we have two orthogonal pairs $a_1,a_2$ and $b_1,b_2$ then
      there are $8$ triple invariants we can compute. In all cases
      we pick $3$ of the flags and order the triple some way and
      then compute.

      \begin{lemma}
        All $8$ triple invariants associated to a pair of orthogonal
        flags
        have the same sign.
      \end{lemma}

      \startproof
      At least for generic choices, we
      can normalize so that one of the pairs is given by
      $$([r,0,1],[-1,0,r]), \hskip 20 pt
      ([-1,0,r],[r,0,1]),$$
      and the other one is given by
      $$([x,y,1],[-x,-y,x^2+y^2])  \hskip 20 pt
      ([-x,-y,x^2+y^2],[x,y,1]))$$
      We compute that
      the triple invariants occur in pairs.  They are
      $t_1$ and $t_2$ and $1/t_1$ and $1/t_2$, where
      $$t_1=\frac{(r-x) \left(r \left(x^2+y^2\right)+x\right)}{(r x+1)
        \left(-r x+x^2+y^2\right)},$$
      $$
t_2=\frac{\left(x^2+y^2+1\right) \left(-r x+x^2+y^2\right) \left(r \left(x^2+y^2\right)+x\right)}{(r-x) (r x+1) \left(\left(x^2+y^2\right)^2+x^2+y^2\right)}.$$
The important thing to notice is that
$$\frac{t_1}{t_2}=\frac{(r-x)^2 \left(x^2+y^2\right)}{\left(r
    x-x^2-y^2\right)^2}>0.$$
Hence $t_1$ and $t_2$ have the same sign.
Taking reciprocals does not change the sign, so
all the triple products have the same sign.
\endproof

Here is the geometric picture. In the negative triple
product case, the lines of the pair $a_1,a_2$ separate the points
of the pair $b_1,b_2$ from each other, and {\it vice versa\/}.
Figure 7.2 shows examples of the positive and the negative
cases.

   \begin{center}
     \resizebox{!}{1.3in}{\includegraphics{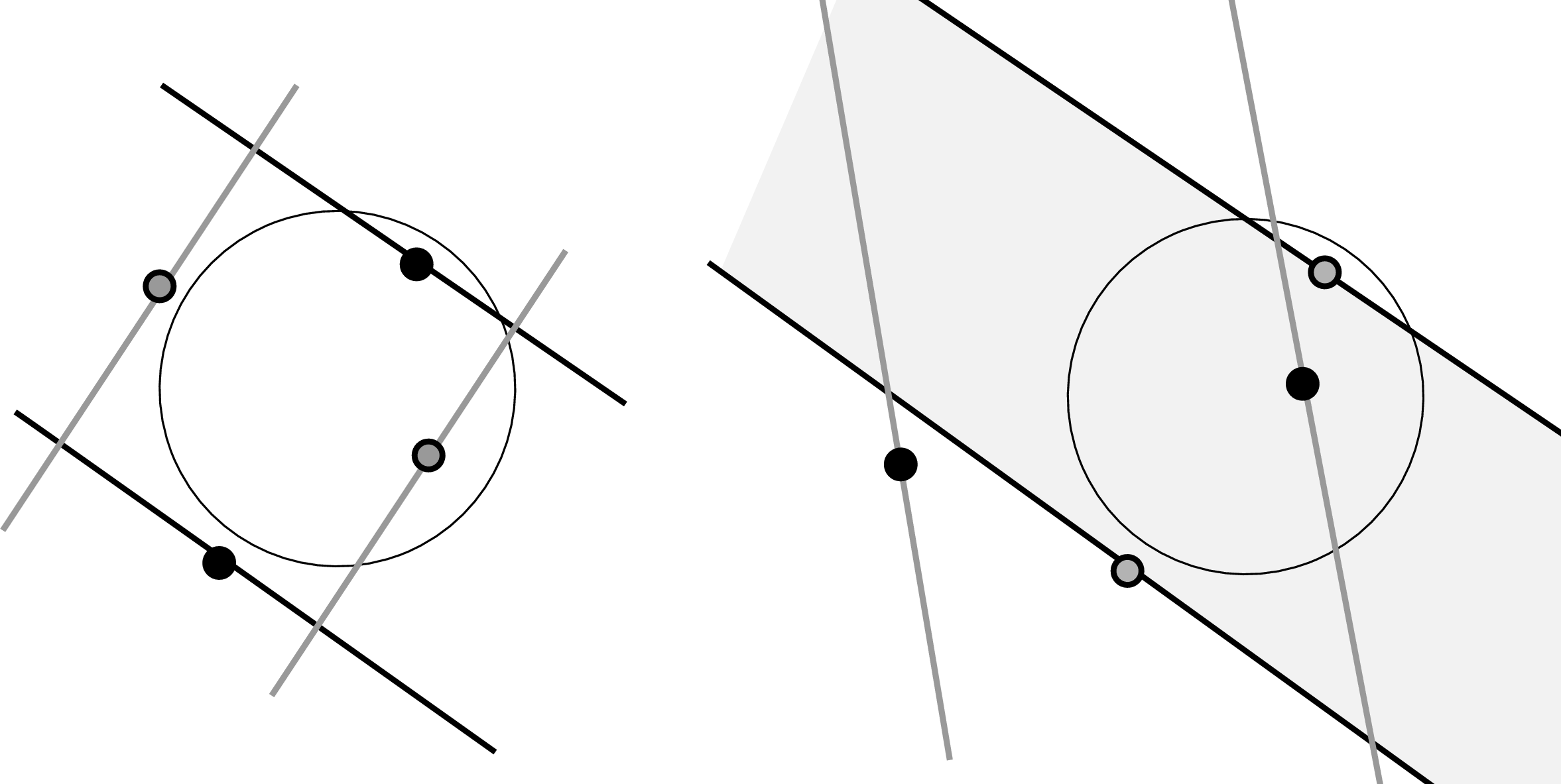}}
    \newline
    {\bf Figure 7.2\/}:   The positive case (left) and
    the negative case (right).
  \end{center}

  More generally, we say that two pairs of flags $(a_1,a_2)$ and
  $(b_1,b_2)$ are {\it separating\/} if the lines of the first pair
  separate the points of the second pair and {\it vice versa\/}.

  \subsection{The Embedding Proof}

  \noindent
  {\bf Transversality:\/}
  Our embedding proof really just uses one familiar
  property of Anosov embeddings, namely
  {\it transversality\/}.  An Anosov representation
  of a group $G$ includes an equivariant map
  \begin{equation}
    \phi: \partial G \to {\cal F\/},
  \end{equation}
  where $\partial G$ is the Gromov boundary
  of $G$ -- in our case a Cantor set -- and
  ${\cal F\/}$ is the flag variety.  The key
  property is that every pair of flags
  of $\phi(\partial G)$ consists of
  transverse flags.   Again, this means that
  the point of one flag does not lie in the
  line of the other. See e.g. [{\bf Lab\/}] or
  [{\bf GW\/}].
  For reference below, we call this collection
  of flags the {\it big collection\/}.
  \newline

  Now we get to the proof.
As we have already mentioned, each prism group
preserves an infinite pattern of flats.  These are
just the orbit of the flats in the initial prism under
the group.  Moreover, each flat in the orbit has
a distinguished geodesic, and so the pattern of
geodesics is embedded provided that the pattern
of flats is embedded.

In [{\bf S0\/}] we proved that the pattern of flats is
embedded when the group is a Pappus modular group.
Our proof in the Anosov case is similar in spirit,
but takes advantage of the transversality property
discussed above.
Corresponding to our prism representation
$\rho$ we have two infinite collections of flags,
one subsuming the other.  We have the
big collection mentioned in the previous
section.   We also have the {\it small collection\/}.
As in the previous chapters we make a choice
of attracting over repelling.  We then take the
attracting flag for the element
$g^2_{\rho}=\rho(\sigma_2\sigma_3\sigma_2\sigma_3)$
and consider its orbit under the group.  This is
the small collection.  The big collection contains the small
collection.

\begin{lemma}
  Every triple of flags in the big collection
  has negative triple invariant.
\end{lemma}

\startproof
This is clearly true for the symmetric Pappus
modular group.  As we move continuously
to other representations,
the invariant cannot change
sign without becoming $0$ along the way.
But if this happens, we have a non-transverse
pair of flags, a contradiction.
\endproof

The flats in our pattern are naturally
associated to the morphed marked boxes
in the orbit.  In the case of Pappus modular
group, the flats are defined in terms of the
tops and the bottoms of the marked boxes.
As we move into the Anosov representations
we define the same kind of association
just by continuity. We call a pair of flags
$(a_1,a_2)$ {\it linked\/} if they are associated
to the same morphed marked box.  This is the
same as saying that they are associated to
the same flat in our pattern.

\begin{lemma}[Separating]
  Let $(a_1,a_2)$ and $(b_1,b_2)$ be
  pairs of linked flags.  Then
  this pair is not separating.
\end{lemma}

\startproof
This also follows from transversality and continuity.
The property is true for shears of the symmetric
Pappus group, as one can see from the nesting of the
geodesics in the hyperbolic plane associated
to these groups.  The general
case follows from continuity.  We we move along
a path of Anosov representation, we can never
acquire the separating property. If we did, we would
encounter a non-transverse pair of flags.
\endproof

Now we will suppose that a pair of flats in our pattern
intersect.   We can move these flats by an isometry
so that their intersection point is the origin.
Now we have a linked pair of orthogonal
flags.  All the triple invariants associated to
these flags must be negative. Hence the linked
pair is separating.  This contradicts the Separating
Lemma. Hence  the flats cannot intersect.
This proves that our pattern of flats is embedded.
This finishes the proof of Theorem \ref{extra}.
To finish the proof of Theorem
\ref{three} we just need to explain the shearing.

\subsection{Shearing}
\label{shear}

Consider two representations $(\Pi,p_1)$ and
$(\Pi,p_2)$ corresponding to points in
the same ray.  The crucial point is that
$p_1$ and $p_2$ lie in the same medial
geodesic $\gamma$ in the same flat $F$ of $\Pi$.
Let $\delta_1=\rho_1(\sigma_2)$ and $\delta_2=\rho_2(\sigma_2)$
denote the elliptic polarities fixing $p_1$ and
$p_2$ respectively.   Let $I$ be the isometry
that translates along $\gamma$, mapping
$p_1$ to $p_2$.    Another crucial point is that
$\delta_k$ conjugates $I$ to $I^{-1}$ because
$\delta_k$ stabilizes $\gamma$ and reverses its
directions. Hence
$$\delta_2=I \circ \delta_1 \circ I^{-1}=I^2 \circ \delta_1.$$
In particular, if $\tau$ is the triangle in the pattern
contained in the prism $\Pi$, then
$$\delta_2(\tau)=I^2(\delta_1(\tau)).$$
Thus the pair $(\tau,\delta_1(\tau))$ is
replaced by the sheared pair
$(\tau,I^2\circ \delta_1(\tau))$.
This is the shearing phenomenon in
Theorems \ref{extra} and \ref{three}.

Finally, the rays are canonically parametrized
as follows.   We say that a point on a given
ray $\gamma$ is $d$ away from the endpoint
if the corresponding representation is obtained
from the Pappus group at the endpoint by
performing a shear of strength $d$.
This pins down the shearing relationship.

It only remains to consider the strengths of
the two kinds of shearing which produce the
same group.  It suffices to consider the generic case.
We revisit the
$r,s,t$ coordinates from \S 5.  If we fix the value of $t$,
then the region $(r,s) \in (0,\infty)^2$ gives us
all the prism groups $(\Pi,p)$ with $p$ varying
in the same flat $F \cong \R^2$ of $\Pi$.

Let $r=\mu s$, as in Equation
\ref{mu}.  The representation corresponding
to $(r,s,t)$ is a Pappus representation.
According to Equation \ref{mu} and
Lemma \ref{ortho2}, the representation
corresponding to $(rd,s/d,t)$ is
a $\sqrt 2 \log d$ shearing of a Pappus representation.
When $d>1$ this is one kind of
shearing and when $d<1$ this is the other.

We compute that eigenvalues
of our element $g^2$ are $-d^2,-d^{-2},1$,
with the flag $(b_1,L_2)$ corresponding
to the $d^2$ eigenvalue.
The eigenvalue set only depends on the shearing strength,
namely $\sqrt 2 \log d$, and not any other
property of the parameters.   The two
prism descriptions of the same
representation give rise to the
same element $g^2$.  Therefore they
correspond to the same strength shears
of both kinds.
This completes the proof of Theorem \ref{three}.

      \section{References}

\noindent
[{\bf Bar\/}]
T. Barbot,
{\it Three dimensional Anosov Flag Manifolds\/},  Geometry $\&$ Topology (2010)
\vskip 12 pt
\noindent
[{\bf BLV\/}], T. Barbot, G. Lee, V. P. Valerio, {\it Pappus's Theorem,
  Schwartz Representations, and Anosov Representations\/},
Ann. Inst. Fourier (Grenoble) {\bf 68\/} (2018) no. 6
\vskip 12 pt
\noindent
[{\bf BCLS\/}] M. Bridgeman, D. Canary, F. Labourie, A. Samburino,
{\it The pressure metric for Anosov representations\/}, Geometry and Functional Analysis {\bf 25\/} (2015)
\vskip 12 pt
\noindent
[{\bf DR\/}] C. Davalo and J. M. Riestenberg, {\it Finite-sided
  Dirichlet domains and Anosov subgroups\/}, arXiv 2402.06408 (2024)
\vskip 12 pt
\noindent
[{\bf FG\/}]. V. Fock and A. Goncharov, {\it Moduli Spaces of local
  systems and
  higher Teichmuller Theory\/}, Publ. IHES {\bf 103\/} (2006)
\vskip 12 pt
\noindent
[{\bf FL\/}] C. Florentino, S. Lawton, {\it The topology of moduli
  spaces of free groups\/}, Math. Annalen {\bf 345\/}, Issue 2 (2009)
\vskip 12 pt
\noindent
[{\bf G1\/}] W. Goldman, {\it Convex real projective structures on compact surfaces\/}, J.
Diff. Geom. {\bf 31\/} (1990)
\vskip 12 pt
\noindent
[{\bf G2\/}], W. Goldman, {\it Mapping Class Group Dynamics on Surface
  Group Representations\/},  arXiv:0509114 (2006)
\vskip 12 pt
\noindent
[{\bf GW\/}] O. Guichard, A. Wienhard, {\it Anosov Representations: Domains of
  Discontinuity and applications\/}, Invent Math {\bf 190\/} (2012)
\vskip 12 pt
\noindent
[{\bf Hit\/}] N. Hitchin, {\it Lie Groups and Teichmuller Space\/}, Topology {\bf 31} (1992)
\vskip 12 pt
\noindent
[{\bf KL\/}] M. Kapovich, B. Leeb, {\it Relativizing characterizations of
  Anosov subgroups, I (with an appendix by Gregory A. Soifer).\/} Groups Geom. Dyn.  {\bf 17\/}  (2023)
  \vskip 12 pt
  \noindent
[{\bf Lab\/}] 
F. Labourie, {\it Anosov Flows, Surface Groups and Curves in Projective Spaces\/},
P.A.M.Q {\bf 3\/} (2007)
\vskip 12 pt
\noindent
[{\bf L\/}] S. Lawton, {\it Generators, relations, and symmeries in
  pairs of $3 \times 3$ unimodular matrices\/}, J. Algebra, 313(2)
(2007)
\vskip 12 pt
\noindent
[{\bf P\/}] R. Penner, {\it The Decorated Teichmuller Theory of
  punctured surfaces\/}, Comm. Math. Pys. {\bf 113\/} (1987)
\vskip 12 pt
\noindent
[{\bf S0\/}] R. E. Schwartz, {\it Pappus's Theorem and the Modular Group\/},
Publ. IHES (1993)
\vskip 12 pt
\noindent
[{\bf S1\/}] R. E. Schwartz, {\it Le Retour de Pappus\/},
KIAS-Springer Lecture Notes (2025) to appear.
See also arXiv 2412.02417
\vskip 12 pt
\noindent
[{\bf S2\/}] R. E. Schwartz, {\it On Pappus and Anosov Representations
  of the Modular Group\/}, preprint (2026)
\vskip 12 pt
\noindent
[{\bf Sil\/}],
  J. Silverman, {\it The arithmetic of dynamical systems\/}, Graduate
  Texts in Mathemtics {\bf 241\/} (2007) Springer
\vskip 12 pt
\noindent
[{\bf T\/}] W. Thurston, {\it The Geometry and Topology of Three
  Manifolds\/}, Princeton University Notes (1978)
\vskip 12 pt
\noindent
[{\bf V\/}], V. P. Valerio, {\it Teorema de Pappus, Representa\c{c}oes de Schwartz e
  Representa\c{c}oes Anosov\/},
Ph. D. Thesis, Federal University of Minas Gerais (2016)
\vskip 12 pt
\noindent
[{\bf W\/}] S. Wolfram et. al., {\it Mathematica\/}, Version 11
Wolfram Research Inc. (2024)

\end{document}